\documentclass[10pt]{amsart}

\usepackage{amssymb, latexsym, amsfonts, pigpen, enumitem, mathrsfs, amsthm, bbm, color}
\usepackage[matrix,arrow,curve]{xy}
\usepackage[plainpages,  urlcolor=blue]{hyperref}
\usepackage{graphicx}   

\newtheorem{theorem}{Theorem}[section]
\newtheorem{thmintro}{Theorem}
\newtheorem*{question}{Question}
\newtheorem{lemma}[theorem]{Lemma}
\newtheorem{proposition}[theorem]{Proposition}
\newtheorem{corollary}[theorem]{Corollary}

\theoremstyle{definition}
\newtheorem{definition}[theorem]{Definition}
\newtheorem{example}[theorem]{Example}
\newtheorem{recollection}[theorem]{Recollection}

\theoremstyle{remark}
\newtheorem{remark}[theorem]{Remark}

\numberwithin{equation}{theorem}

\makeatletter
\def\blfootnote{\gdef\@thefnmark{}\@footnotetext}
\makeatother

\newcommand{\GW}{{\mathrm{GW}}}

\newcommand{\Sm}{{\mathrm{Sm}}}
\newcommand{\I}{{\mathrm{I}}}
\renewcommand{\H}{\mathrm{H}}

\newcommand{\W}{\mathrm{W}}
\newcommand{\K}{\mathrm{K}}
\newcommand{\M}{\mathrm{M}}
\newcommand{\MW}{\mathrm{MW}}
\newcommand{\KM}{\mathrm{K}^{\M}}
\newcommand{\kM}{\mathrm{k}^{\M}}
\newcommand{\KMW}{\mathrm{K}^{\MW}}
\newcommand{\CH}{\mathrm{CH}}
\newcommand{\Ch}{\mathrm{Ch}}

\newcommand{\A}{{\mathbb{A}}}
\newcommand{\C}{{\mathbb{C}}}
\newcommand{\R}{{\mathbb{R}}}

\newcommand{\Z}{{\mathbb{Z}}}

\newcommand{\PP}{{\mathbb{P}}}

\newcommand{\QQ}{{\mathbb{Q}}}
\newcommand{\rank}{\operatorname{rank}}

\newcommand{\struct}{\mathcal{O}}

\newcommand{\GL}{{\mathrm{GL}}}

\newcommand{\Um}{{\mathrm{Um}}}
\newcommand{\diag}{{\mathrm{diag}}}

\newcommand{\Spec}{\operatorname{Spec}}
\newcommand{\tr}{\operatorname{tr}}
\newcommand{\Zar}{{\mathrm{Zar}}}

\newcommand{\id}{\mathrm{id}}
\newcommand{\Id}{\mathrm{Id}}

\newcommand{\chark}{\operatorname{char}}

\newcommand{\CHW}{{\widetilde{\mathrm{CH}}}{}}

\newcommand{\coker}{\operatorname{coker}}

\newcommand{\Gr}{{\mathrm{Gr}}}

\newcommand{\T}{{\mathfrak{T}}}

\newcommand{\Vect}{\mc{V}}
\newcommand{\image}{\operatorname{Im}}

\newcommand{\lmod}{\textrm{-}\mathrm{mod}}
\newcommand{\Coh}{\operatorname{Coh}}
\newcommand{\Nrd}{\operatorname{Nrd}}
\newcommand{\Sym}{\operatorname{S}}

\newcommand{\mc}[1]{\mathcal{#1}}
\newcommand{\mr}[1]{\mathrm{#1}}

\newcommand{\Lodd}[1]{\Lambda^{\mr{odd}}({#1})}
\newcommand{\Lev}[1]{\Lambda^{\mr{ev}}({#1})}
\newcommand{\Sus}[1]{\mr{S}({#1})}

\begin{document}

\title{The vector fields problem on smooth affine quadrics}

\author{Alexey Ananyevskiy}
\address{Mathematisches Institut der Universit\"at M\"unchen, Theresienstr. 39, D-80333 M\"unchen, Germany}
\email{alseang@gmail.com}

\author{Jean Fasel}
\address{Univ. Grenoble Alpes, CNRS, Institut Fourier, 38000 Grenoble, France}
\email{jean.fasel@univ-grenoble-alpes.fr}

\date{}

\begin{abstract}
We study the problem of trivializing the tangent bundle to smooth affine quadrics over a perfect field $k$ of characteristic different from $2$. We obtain complete answers and explicit formulas in dimensions $2$ and $3$, and necessary and sufficient conditions in higher dimensions. In particular, we show that in the even-dimensional case the tangent bundle is trivial if and only if it admits a nowhere vanishing section, which in some cases may be detected by the group of values of the quadratic form defining the quadric but is in general more involved as evidenced by the example of an anisotropic Albert quadric that we study extensively. These general criteria are complemented with a construction of explicit families of affine quadrics with trivial tangent bundle. The proofs are based on computations with cohomological invariants valued in Milnor(--Witt) $\K$-cohomology groups of quadrics and an analysis of the latter groups using in particular a theorem of Rost that we reprove in an appendix.
\end{abstract}

\maketitle

\tableofcontents

\section{Introduction}

In topology, the study of the tangent bundles of spheres has a rich history.  The famous hairy ball theorem states that there is no nowhere-vanishing (continuous) section of the tangent bundle of an even dimensional sphere. The situation is more complicated for spheres of odd dimension. Using the fact that they are unit spheres in the division algebras $\C$, $\mathbb{H}$ and $\mathbb{O}$, it is straightforward to see that $S^1,S^3$ and $S^7$ are parallelizable (i.e. their tangent bundle is trivial). Bott-Milnor \cite{BM58} and Kervaire \cite{Ker58} proved independently that these are the only parallelizable spheres. Nevertheless, it is natural to ask for the maximal number of linearly independent vector fields on $S^{n-1}$. If $n=(2a+1)2^{c+4d}$ with $0\leq c\leq 3$, set $\rho(n)=2^c+8d$. Using the Hurwitz-Radon matrices, Eckmann \cite{Eck42} and Hopf \cite{Hop40} independently observed that $S^{n-1}$ carries $\rho(n)-1$ independent vector fields, while Adams proved later on in \cite{Ada62}, using the eponymous operations, that this number was indeed maximal.

From its inception, the theory of projective modules has been meant as the algebraic version of vector bundles, and it is not surprising that the above questions attracted some attention in the algebraic world. If 
\[
A_n:=\R[x_1,\ldots,x_n]/(\sum x_i^2-1)
\] 
denotes the coordinate ring of the real algebraic sphere $S^{n-1}_{\R}$, it is of course possible (following \cite{Swa62}) to extend scalars to the ring $C(S^{n-1},\R)$ of continuous functions on $S^{n-1}$, and use the topological results to obtain the same answers, both positive and negative. For instance, the tangent bundle to $S^{n-1}_{\R}$ has no nowhere vanishing global sections if $n-1$ is even. In the same vein, the construction of the independent vector fields being algebraic, there are $\rho(n)-1$ linearly independent vector fields on $S^{n-1}_{\R}$. Although straightforward, the proofs of these results are a bit unsatisfactory for an algebraically minded mathematician and the question to find algebraic proofs was raised early on. One example is provided by M. Kong in \cite{Kon77}, where he was able to define the Euler class of a projective module endowed with an inner product using (algebraic) de Rham cohomology, and give an algebraic proof that the tangent bundle to $S^{n-1}_{\R}$ has no nowhere vanishing section for $n$ odd.

More recently, Bhatwadekar--Sridharan developed, following an original idea of Nori, a purely algebraic theory of Euler classes hosted in the so-called Euler class groups \cite{BS98}, \cite{BS99}. Around the same time, Barge--Morel defined the Chow--Witt groups \cite{BM00} (see \cite{Fas08} for a full development) and an Euler class with coefficients in these groups. In both theories, the triviality of the Euler class of a top rank vector bundle over a smooth affine variety is a necessary and sufficient condition to have a nowhere vanishing section. The Chow--Witt groups are now understood as part of Milnor--Witt motivic cohomology, one of the main candidates for an algebraic analogue of singular cohomology. In particular, many of the basic structural results of the latter are also available in the former and it is possible to translate these cohomological techniques to better understand algebraic vector bundles. 

As is customary in mathematics, these recent developments renewed interest in tangent bundles to algebraic spheres, taken in a broad sense. Indeed, if $k$ is a field, the split affine quadrics are particular cases of motivic spheres \cite{ADF17} and non singular affine quadrics over $k$ are forms of motivic spheres, in the sense that they become isomorphic to a motivic sphere after a suitable extension of the base field. In this context, the following result of Levine and the first author of this article can be seen as a first step towards the study of the vector fields problem on smooth affine quadrics. 

\begin{thmintro}[{\cite[Theorem~4.1]{AL25}}]
    Let $q=\sum_{i=1}^{n+1}a_ix_i^2$ be a regular quadratic form over a perfect field $k$ of characteristic different from $2$, and let $Q^o$ be the affine quadric given by the equation $\sum_{i=1}^{n+1}a_ix_i^2=1$.
	Then the following statements hold.
	\begin{enumerate}
		\item 
		If $n$ is odd then the tangent bundle $T_{Q^o}$ has a nowhere vanishing section.
		\item 
		If $n>0$ is even and the tangent bundle $T_{Q^o}$ has a nowhere vanishing section then 
		\[-\prod_{i=1}^{n+1} a_i \in [D(q)^2].\]
		\item 
		If $n>0$ is even and $Q^o(k)\neq \emptyset$ then the tangent bundle $T_{Q^o}$ has a nowhere vanishing section if and only if $-1 \in [D(q)]$.
	\end{enumerate}
\end{thmintro}

This result can be thought of as an algebraic analogue of the hairy ball theorem, and the purpose of this article is to push this result further towards its topological analogue. In the paper we work over a perfect field $k$ of characteristic $\neq 2$. Let us start with discussing the case of even dimensional quadrics. The first result we want to emphasize is the surprising fact that the tangent bundle has a nowhere vanishing section if and only if it is trivial  (Theorem \ref{thm:even_parallel} below).

\begin{thmintro}
    Let $n>0$ be an even integer, and let $Q^o$ be the smooth affine quadric given by the equation $\sum_{i=1}^{n+1}a_ix_i^2=1$. The following conditions are equivalent.
    \begin{enumerate}
        \item $T_{Q^o}$ has a nowhere vanishing section,
        \item$T_{Q^o}$ is trivial.
    \end{enumerate}
    Let $Q$ be the projective quadric given by the equation $\sum_{i=1}^{n+1}a_ix_i^2=x_0^2$, if the pushforward homomorphism $\pi_*\colon \CHW_0(Q)\to \GW(k)$ is injective, then the above are also equivalent to
    \begin{enumerate}[resume]
        \item $-a_1 \in [D(q)^2]$.
    \end{enumerate}
\end{thmintro}
\noindent
The condition on the pushforward $\pi_*\colon \CHW_0(Q)\to \GW(k)$ ensures that one can detect the triviality of the Euler class of the tangent bundle using the Euler characteristic. This condition is satisfied for instance if $Q$ is isotropic (in which case the condition $(3)$ boils down to $-1 \in [D(q)]$), but fails in general. We spend a considerable amount of effort in studying this map for any quadric (regardless of the parity of its dimension), and prove in Theorem \ref{thm:pushforwardequivalence} that its kernel is isomorphic to the kernel of the pushforward homomorphism $\pi_*^{\mathrm{W}}\colon \H^n(Q,\mathrm{W}\otimes\omega_{Q/k})\to \mathrm{W}(k)$ which is slightly simpler to study. A direct consequence is that $\pi_*\colon \CHW_0(Q)\to \GW(k)$ is injective for $1\le \mathrm{dim}(Q)\leq 3$, but we provide a counter-example in dimension $4$ under the guise of an anisotropic Albert quadric. In relation with the above theorem, we also study the Euler class of its tangent bundle, and prove in Section \ref{sec:SuslinEuler} that it is nontrivial even though the given quadric satisfies condition (3). Consequently \cite[Theorem~4.1.(2)]{AL25} gives only a necessary condition for the existence of a nowhere vanishing section, but in general not a sufficient one. To conclude this discussion on even dimensional quadrics, let us just say that we also provide explicit formulas for the trivialization of the tangent bundle in dimension $2$ (Theorem~\ref{thm:hedge2} and Examples~\ref{ex:dim2_1}, \ref{ex:dim2_2}, and~\ref{ex:dim2_3}).

The first obvious thing to say in case of odd dimensional smooth affine quadrics is the well-known fact that the tangent bundle has a nowhere vanishing section. Regarding the triviality of the tangent bundle, we offer the following results gathered from Theorems \ref{thm:n_odd_necessary} and \ref{thm:n_odd_sufficient}.

\begin{thmintro}
Let $n\geq 3$ be an odd integer, and let $Q^o$ be the smooth affine quadric given by the equation $\sum_{i=1}^{n+1}a_ix_i^2=1$, and $Q$ be the projective quadric given by the equation $\sum_{i=1}^{n+1}a_ix_i^2=x_0^2$. The following statements hold:
\begin{enumerate}
\item If $(-1)^{\frac{n+1}{2}}\prod_{i=1}^{n+1} a_i \notin (k^\times)^2$ and  $T_{Q^o}$ is trivial, then $-a_1\in [D(q)^2]$.
\item If the pushforward $\pi_*\colon \H^n(Q,\K^\M_{n+1})\to k^\times$ is injective and $-a_1\in [D(q)^2]$, then the tangent bundle $T_{Q^o}$ is trivial.
\end{enumerate}
\end{thmintro}   
\noindent
Again, if the quadric $Q$ is isotropic, then the pushforward is injective and the second condition boils down to $-1\in [D(q)]$. The injectivity of $\pi_*$ for arbitrary quadrics was studied by Rost and Chernousov--Merkurjev, see \cite[Proposition~2]{Ros88} and \cite[Corollary~6.3]{CM01}. Combining the above theorem with the explicit formulas from Example~\ref{ex:dim3} we obtain a complete answer in dimension $3$, where we prove that $T_{Q^o}$ is trivial if and only if either $\prod_{i=1}^{4} a_i\in (k^\times)^2$ or $-a_1\in [D(q)^2]$ (Corollary \ref{cor:dim3}).

Finally, using the stable birational geometry of quadrics we construct the following families of examples of smooth affine quadrics with trivial tangent bundles (Theorem~\ref{thm:vanishing_birational}). 
\begin{thmintro}
    Let $Q^o$ be the smooth affine quadric given by the equation $\sum_{i=1}^{n+1}a_ix_i^2=1$, let $\tilde{q}:=-x_0^2+q$ and suppose that $\tilde{q}$ is anisotropic and the first higher Witt index satisfies $i_1(\tilde{q})\ge 2$. Then $T_{Q^o}$ is trivial.
\end{thmintro}
\noindent
The examples of such $\tilde{q}$ include non-minimal Pfister neighbors \cite[\S 3, Corollary 2]{Hof95}.

To prove our results, we use cohomological methods under the form of computations with the groups of (co-)cycles with coefficients in Milnor--Witt cycle modules (in the sense of \cite{Fel20}), most notably Milnor(--Witt) $\K$-theory. A nontrivial ingredient is provided by a special case of a theorem due to Rost, describing the group $\H^{n-1}(Q,\kM_n)$ for a smooth projective quadric of dimension $n$ in terms of pushforward maps of much simpler groups. Since his result remains unpublished, we provide a proof of the required statement in Appendix \ref{sec:Rosttheorem}, suggested by an LLM (see our statement below).  One consequence of the method of the proof is a particularly nice description in Lemma \ref{lem:baby_rost} of the (complicated) group of $1$-cycles on a smooth quadric $Q$ by providing easy geometric origins for the generators. We counterbalance these rather abstract methods by exhibiting in Appendix~\ref{sec:explicitformulas} explicit formulas for the trivialization of the tangent bundle whenever we were able to do it.

We finish the introduction with the following question. 
\begin{question}[vector fields problem for smooth affine quadrics] Let $Q^o$ be a smooth affine quadric of dimension $n$. Can we describe the maximal number $m$ of linearly independent sections of the tangent bundle $T_{Q^o}$ in terms of classical invariants attached to the underlying quadratic form?
\end{question}
\noindent Currently we have a complete answer in dimensions $\le 3$, we know that in the even dimension $m=0$ or $m=n$ which in favorable cases can be detected using the group of values, while in the odd dimension we classically know that $m\ge 1$ and have some necessary and sufficient conditions for $m=n$. There are also results for certain families of quadrics and over special classes of field, but in general the question is wildly open.

\subsection*{Notations}

Throughout the paper we employ the following notation.
\\[5pt]
	\begin{tabular}{l|l}
	$k$ & a perfect field with $\chark k\neq 2$ \\
    $q$ & a regular quadratic form over $k$ \\
    $Q^o$ & the affine quadric $\{q=1\}$\\
    $R$, $R(q)$ & the ring of regular functions $k[Q^o]$ on $Q^o$ \\
    $\Um_n(R)$ & the set of length $n$ unimodular rows over $R$\\
	$\langle a_1,a_2,\dots,a_n\rangle$ & the quadratic form $a_1x_1^2+a_2x_2^2+\dots + a_nx_n^2$ \\
    $i_1(q)$ & the first higher Witt index\\
    $D(q)$, $D(q)^2$ & the sets of non-zero values and products of pairs of non-zero values of $q$\\
    $[D(q)]$, $[D(q)^2]$ & the subgroups of $k^\times$ generated by the respective subsets\\
    $\Delta$, $\Delta(q)$ & signed discriminant $(-1)^{\frac{n(n-1)}{2}}\prod_{i=1}^n a_i$ for a quadratic form $q=\sum_{i=1}^n a_ix_i^2$\\    
    $C(q)$, $C_0(q)$ & Clifford algebra of $q$ and its even part\\
    $N_{L/k}$& norm map for a finite field extension $L/k$\\
    $\Nrd_{A/k}$ & reduced norm map for a central simple algebra $A/k$\\
    $\Sm_k$ & the category of smooth varieties over $k$ \\
    $\CH^*$, $\Ch^*$ & Chow ring and Chow ring modulo $2$\\ 
    $\KMW_*$, $\KM_*$, $\kM_*$ & Milnor--Witt $\K$-theory, Milnor $\K$-theory and Milnor $\K$-theory modulo $2$ \\
    $\W$, $\I^*$, $\overline{\I}{}^*$ & Witt ring, powers of the fundamental ideal in $\W$ and $\overline{\I}{}^*:=\I^*/\I^{*+1}$ \\
	\end{tabular}
\\

\subsection*{Acknowledgments}    
The authors are grateful to Heng Xie for useful conversations on Hermitian $\K$-theory of quadrics, to Stephen Scully for pointing out that Albert forms usually have nontrivial Witt-kernels, and to Alexander Merkurjev for the discussions about explicit formulas and cycles on quadrics. The work of the first author is supported by DFG Heisenberg grant AN 1545/1-1 and by the DFG research grant AN 1545/4-1. Part of the work was done while the author was in residence at the Simons Laufer Mathematical Sciences Institute in Berkeley, California, during the Fall 2026 semester, where he was partially funded by the National Science Foundation under Grant DMS-2424139. The second author wishes to thank the CAS and the LMU in Munich for supporting a visit in 2025 during which part of the article was discussed. 

\subsection*{Statement about the use of AI} LLMs were used to find explicit formulas in dimension~$2$ and to search the literature. While we were looking for a published account of Rost's result mentioned in \cite[Remark~4.6]{Mer00}, an LLM suggested a sketch of a proof of a particular case needed in the paper. We expanded this sketch and included the proof as Theorem~\ref{thm:rost}. LLMs were also used to check the manuscript for mistakes and consistency of notation. The authors take full responsibility for the final work.

\section{Preliminaries on quadratic forms and unimodular rows}
In this section we settle notation related to quadratic forms, unimodular rows and cohomology groups with coefficients in homotopy modules, and recall some basic properties thereof.

\begin{definition}
	Let $q$ be a regular quadratic form over $k$. We say that $q$ is \textit{isotropic} over $k$ if $q=0$ has a non-trivial solution over $k$, or, equivalently, if the projective quadric $\{q=0\}$ has a rational point. Otherwise $q$ is \textit{anisotropic}. For a field extension $K/k$ we denote by $q_K$ the quadratic form $q$ extended to $K$.
    We use the following notation:
	\begin{itemize}
		\item 
		$D(q)\subseteq k^\times$ is the set of nonzero values of $q$,
		\item
		$D(q)^2=\{a\cdot b\mid a,b\in D(q)\}$ is the set of products of pairs of nonzero values of $q$,
		\item
		$[D(q)]$, $[D(q)^2]$ are the multiplicative subgroups of $k^\times$ generated by the respective sets. 
	\end{itemize}
\end{definition}

The group $[D(q)^2]$ admits an alternative description as the group of norms associated to $q$:
\begin{lemma}[{\cite[Lemma~2.2]{CTS93}}] \label{lem:CTS}
    Suppose $\dim q\ge 2$. Then $[D(q)^2]$ coincides with the subgroup of $k^\times$ generated by the norms $N_{L/k}(a)$, $a\in L^\times$, for $L/k$ finite such that $q_L$ is isotropic.
\end{lemma}

\begin{definition}
	For $a_1,a_2,\dots,a_n\in k^\times$ we denote by 
    \begin{itemize}
        \item  $\langle a_1,a_2,\dots, a_n \rangle:=\sum_{i=1}^n a_ix_i^2$ the associated diagonal quadratic form of dimension $n$,
        \item $\langle\langle a_1,a_2,\dots, a_n \rangle\rangle:=\prod_{i=1}^n \langle 1, -a_i\rangle$ the associated \textit{$n$-fold Pfister form} of dimension $2^n$.
    \end{itemize}
    A regular quadratic form $q$ over $k$ is called a \textit{Pfister neighbor} if there exists $a\in k^\times$ such that $\langle a\rangle \cdot q$ is a subform of an $n$-fold Pfister form with $2^{n-1}<\dim q$. Note that the Pfister form containing $\langle a\rangle \cdot q$ for a Pfister neighbor $q$ is unique \cite[Proposition~X.4.17]{Lam05}. For a Pfister neighbor $q$ with associated Pfister form $\phi$ one has $[D(q)^2]=D(\phi)$, see e.g. \cite[Lemma~3.15]{AL25}.
\end{definition}

\begin{definition}
    For a regular quadratic form $q$ over $k$ we denote by $C(q)$ and $C_0(q)$ the associated Clifford algebra and its even part. Recall that if $\dim q$ is even then $C(q)$ is a central simple algebra over $k$ and if $\dim q$ is odd then $C_0(q)$ is a central simple algebra over $k$.
\end{definition}

\begin{definition}
    Let $R$ be a commutative ring. We say that a row $u=(u_1,u_2,\hdots, u_n)\in R^n$ is
    \begin{itemize}
        \item \textit{unimodular} if $Ru_1+Ru_2+\hdots + Ru_n=R$,
        \item \textit{completable over $R$} if $u$ may be completed to a basis of $R^n$, or, equivalently, if there exists $g\in \GL_n(R)$ with $(u_1,u_2,\hdots, u_n)$ being the first row.
    \end{itemize}
    We denote by $\Um_n(R)$ the set of length $n$ unimodular rows over $R$. Multiplication on the right endows $\Um_n(R)$ with the action of $\GL_n(R)$. A row is completable if it lies in the orbit $(1,0,\hdots,0)\cdot \GL_n(R)$, in particular, a completable row is unimodular. A unimodular row in general may fail to be completable.
\end{definition}

\begin{definition} \label{def:tangent-normal}
    For a regular quadratic form $q=\sum_{i=1}^{n+1}a_ix_i^2$ over $k$ we denote by
    \[
    R(q):=k[x_1,x_2,\hdots,x_{n+1}]/(q-1),\quad Q^o:=\Spec R(q)
    \]
    the associated affine quadric and its ring of regular functions. We usually write $R:=R(q)$ when the quadratic form is clear from the context.
    
    The embedding $Q^o\subseteq \mathbb{A}^{n+1}$ gives rise to a short exact sequence of vector bundles
    \[
    0\to T_{Q^o} \to \struct^{\oplus n+1}_{Q^o} \to \mc{N} \to 0
    \]
    with $\mc{N}$ being the normal bundle of the embedding. The polynomial $q$ gives a trivialization of~$\mc{N}$, so the above sequence of vector bundles may be rewritten as a short exact sequence
    \[
    0\to \T\to R^{ n+1} \xrightarrow{\nu} R \to 0
    \]
    of projective modules over $R$, where the epimorphism $\nu$ is given by
    \[
    \nu(u_1,u_2,\hdots,u_{n+1})=a_1x_1u_1+a_2x_2u_2+\hdots +a_{n+1}x_{n+1}u_{n+1}.
    \]
    The sequence splits via 
    \[
    s\colon R\to R^{ n+1}, \quad s(u)=(x_1u,x_2u,\hdots,x_{n+1}u),
    \]
    yielding the decompositions
    \begin{equation} \label{eq:tangent_stably_trivial}
    R^{ n+1}= \T\oplus R\cdot(x_1,x_2,\hdots,x_{n+1}), \qquad      \struct^{\oplus n+1}_{Q^o} \cong T_{Q^o} \oplus \mc{N}.
    \end{equation} 
    The tangent bundle $T_{Q^o}$ has a nowhere vanishing section if and only if there exists a unimodular row $u=(u_1,u_2,\hdots,u_{n+1})\in R^{n+1}$ such that $\nu(u)=0$.
\end{definition} 

\begin{lemma} \label{lem:completable_trivial}
    Let $q=\sum_{i=1}^{n+1}a_ix_i^2$ be a regular quadratic form, $R:=R(q)$ and $Q^o:=\Spec R$. The following are equivalent.
    \begin{enumerate}
        \item $(x_1,x_2,\hdots,x_{n+1})\in R^{n+1}$ is completable over $R$,
        \item the tangent bundle $T_{Q^o}$ is trivial.
    \end{enumerate}
\end{lemma}
\begin{proof}
    Completability of $(x_1,x_2,\hdots,x_{n+1})$ means that $(x_1,x_2,\hdots,x_{n+1})$ can be completed to a basis of $R^{ n+1}$, which is equivalent to the triviality of $T_{Q^o}$ by the decomposition~\eqref{eq:tangent_stably_trivial}.
\end{proof}

\section{Dimension $2$}\label{sec:Dim2}

In this section we establish a necessary and sufficient condition for an affine quadric of dimension $2$ to admit a nowhere vanishing vector field, or, equivalently, to have a trivial tangent bundle. This is done by combining the necessary condition of \cite[Theorem~4.1]{AL25} with a construction of an explicit trivialization based on computations in the associated quaternion algebra.

\begin{definition}
    A \textit{quaternion algebra} over a field $k$ is a central simple algebra over $k$ of degree $2$. It can be shown that every quaternion algebra $\mr{H}$ admits a \textit{quaternion basis} $\{1,i,j,ij\}$ subject to the relations $i^2,j^2\in k^\times$, $ij=-ji$. If $i^2=a$ and $j^2=b$ the quaternion algebra $\mr{H}$ is denoted
    \[
    \mr{H}=(a,b)_k.
    \]
    The algebra $\mr{H}$ has a canonical normalized reduced trace homomorphism $\operatorname{Trd}_{\mr{H}/k}\colon \mr{H}\to k$ and a reduced norm quadratic map $\operatorname{Nrd}_{\mr{H}/k}\colon \mr{H}\to k$. In the quaternion basis $\{1,i,j,ij\}$ we have
    \[
    \operatorname{Trd}_{\mr{H}/k} (x_0+x_1i+x_2j+x_3ij) = x_0,\quad \Nrd_{\mr{H}/k} (x_0+x_1i+x_2j+x_3ij) = x_0^2-ax_1^2-bx_2^2+abx_3^2.
    \]
    We denote
    \[
    \mr{H}^0:=\{h\in \mr{H}\,|\, \operatorname{Trd}_{\mr{H}/k} h =0\}.
    \]
    For $h\in \mr{H}^0$ one has $h^2=-\Nrd_{\mr{H}/k} h \in k$. If $u,v \in \mr{H}^0$ are invertible and satisfy $uv \in \mr{H}^0$ then $\{1,u,v,uv\}$ form a quaternion basis.
    
    A quaternion algebra is \textit{split} if it is isomorphic to $\mr{M}_2(k)$, in this case $\operatorname{Trd}_{\mr{H}/k} g = \frac{1}{2}\tr g$ and $\Nrd_{\mr{H}/k} g = \det g$. A non-split quaternion algebra is necessarily a division algebra, and every non-split quaternion algebra splits after a quadratic field extension.
\end{definition}

\begin{lemma} \label{lem:quaternion_basis}
    Let $\mr{H}=(a,b)_k$ be a non-split quaternion algebra, $\Delta\in k^\times$ and suppose that there exists $h\in \mr{H}$ such that $\Nrd_{\mr{H}/k} h =\Delta$. Then there exists a quaternion basis $\{1,u,v,uv\}$ and $r\in k^\times$ such that $\Nrd_{\mr{H}/k} (r+uv) = \Delta$.
\end{lemma}
\begin{proof}
    Since $\mr{H}$ is not split, $\mr{H}$ is a division algebra and $k$ is infinite. The smooth affine quadric $X:=\{x_0^2-ax_1^2-bx_2^2+abx_3^2=\Delta\}$ has a rational point by assumption, thus it is a rational variety. Since $k$ is infinite, the set of rational points on a rational variety is dense in the Zariski topology. It follows that there exists $t\in \mr{H}\setminus \mr{H}^0$ such that $\Nrd_{\mr{H}/k}t=\Delta$ and $\operatorname{Trd}_{\mr{H}/k}t\neq t$. Let $w:= t - \operatorname{Trd}_{\mr{H}/k}t$ and note that $w\in \mr{H}^0\setminus\{0\}$. Choose some $u\in \mr{H}^0\setminus\{0\}$ such that $uw\in \mr{H}^0$, and put $v:=u^{-1}w$, $r:=\operatorname{Trd}_{\mr{H}/k}t$. The claim follows.
\end{proof}

\begin{remark}
    The choice of the quaternion $t$ in the proof of Lemma~\ref{lem:quaternion_basis} may be made explicit using the stereographic projection from the point $h\in X$ and choosing an appropriate rational point in the image.
\end{remark}

\begin{theorem} \label{thm:hedge2}
    Let $q=a_1x^2+a_2y^2+a_3z^2\in k[x,y,z]$ for some $a_1,a_2,a_3\in k^\times$, and denote
    \[
    a:=-a_2a_3,\quad b:=-a_1a_3,\quad \Delta:=-a_1a_2a_3,\quad R:=k[x,y,z]/(q-1),\quad Q^o:=\Spec R.
    \]
    Then the following are equivalent:
    \begin{enumerate}
        \item $\Delta\in [D(q)^2]$,
        \item $x_0^2-ax_1^2-bx_2^2+abx_3^2=\Delta$ has a solution over $k$,
        \item the row $(x,y,z)$ is completable over $R$,
        \item the tangent bundle $T_{Q^o}$ is trivial,
        \item the tangent bundle $T_{Q^o}$ has a nowhere vanishing section.
    \end{enumerate}
\end{theorem}
\begin{proof}
    $(1) \Leftrightarrow (2)$. The form $a_3\langle a_1,a_2,a_3\rangle \simeq\langle -b, -a,1\rangle$ is a subform of the Pfister form $\langle\langle a,b\rangle\rangle$, thus $[D(q)^2]=D(x_0^2-ax_1^2-bx_2^2+abx_3^2)$ by \cite[Lemma~3.15]{AL25}.

    $(3) \Leftrightarrow (4)$. This is Lemma~\ref{lem:completable_trivial}.

    $(4) \Leftrightarrow (5)$. The decomposition~\eqref{eq:tangent_stably_trivial} yields an isomorphism $\struct_{Q^o} \cong \Lambda^2 T_{Q^o}$. If $T_{Q^o}$ has a nowhere vanishing section then $T_{Q^o}\cong \struct_{Q^o}\oplus L$ for some line bundle, and it follows that 
    \[
    \struct_{Q^o} \cong \Lambda^2 T_{Q^o}\cong \Lambda^2(\struct_{Q^o}\oplus L) \cong L.
    \]
    The other direction is clear.

    $(5) \Rightarrow (1)$. Follows from \cite[Theorem~4.1.(2)]{AL25}.

    $(2) \Rightarrow (3)$.  If $q$ is isotropic then the claim follows from \cite[2.1]{Rev79}, see also Example~\ref{ex:Revoy} below for the explicit formulas. Thus we may assume $q$ to be anisotropic. Consider the algebra
    \[
    \tilde{R}:=k[\tilde{x},\tilde{y},\tilde{z}]/(-a \tilde{x}^2-b \tilde{y}^2 +ab \tilde{z}^2+\Delta)
    \]
    and note that there is an isomorphism $R\xrightarrow{\simeq} \tilde{R}$ given by $(x,y,z)\mapsto (a_1^{-1}\tilde{x}, a_2^{-1}\tilde{y},\tilde{z})$, in particular, $(x,y,z)$ is completable over $R$ if and only if $(\tilde{x},\tilde{y},\tilde{z})$ is completable over $\tilde{R}$. Consider the quaternion algebra $\mr{H}:=(a,b)_k$ and denote
    \[
        \mr{H}^0(\tilde{R}):= \tilde{R} \otimes_k \mr{H}^0 = \tilde{R}i\oplus \tilde{R}j\oplus \tilde{R}ij.
    \]
    Completability of $(\tilde{x},\tilde{y},\tilde{z})$ over $\tilde{R}$ means that $\tilde{x}i+\tilde{y}j+\tilde{z}ij$ can be completed to a basis of the $\tilde{R}$-module $\mr{H}^0(\tilde{R})$. Since $q$ is anisotropic, $\mr{H}$ is non-split.
    Lemma~\ref{lem:quaternion_basis} yields that there exists a quaternion basis $\{1,u,v,uv\}$ of $\mr{H}$ and $r\in k^\times$ such that $\Nrd_{\mr{H}/k} (r+uv) = \Delta$. Note that we have
    \[
        \mr{H}^0(\tilde{R})= \tilde{R}u\oplus \tilde{R}v\oplus \tilde{R}uv.
    \]    
    Changing the basis $\{i,j,ij\}$ to $\{u,v,uv\}$ we obtain
    \[
        \tilde{x}i+\tilde{y}j+\tilde{z}ij = \hat{x}u+\hat{y}v+\hat{z}uv,
    \]
    where $(\hat{x},\hat{y},\hat{z}) = (\tilde{x},\tilde{y},\tilde{z}) (S^T)^{-1}$ for the matrix $S\in \GL_3(k)$ satisfying $(u,v,uv)=(i,j,ij)S$. Let 
    \[
        \alpha:= u^2,\quad \beta:= v^2,
    \]
    and note that $\hat{x}$, $\hat{y}$ and $\hat{z}$ satisfy the relation $-\alpha \hat{x}^2-\beta \hat{y}^2 +\alpha \beta \hat{z}^2 = -\Delta$. Consider the matrix
     \[
        M=
        \begin{pmatrix}
            \hat{x} & \hat{y} & \hat{z} \\
            1-\frac{\alpha\beta}{\Delta} \hat{z} & -\frac{\alpha r}{\Delta} \hat{z} & -\frac{2\alpha}{\Delta}\hat{x} + \frac{\alpha \beta - r^2}{r \Delta}\hat{y}\\
            \frac{\beta r}{\Delta} \hat{z} & 1-\frac{\alpha\beta}{\Delta} \hat{z} & -\frac{\alpha \beta - r^2}{r \Delta } \hat{x} -\frac{2\beta}{\Delta} \hat{y}
        \end{pmatrix}.
    \]
    A straightforward computation shows that
    \[
        \det M = 2+ \frac{2-\hat{z}}{\Delta}\left( \alpha \hat{x}^2+\beta \hat{y}^2-\alpha\beta \hat{z}^2-\Delta \right) = 2.
    \]
    Then $M\cdot(u,v,uv)^T=(\hat{x}u+\hat{y}v+\hat{z}uv,\hdots,\hdots)^T$ is a basis of $\mr{H}^0(\tilde{R})$ over $\tilde{R}$, so $\tilde{x}i+\tilde{y}j+\tilde{z}ij$ is completable to a basis, and hence $(x,y,z)$ is completable over $R$. An explicit matrix with the first row being $(x,y,z)$ is given by $M\cdot S^T\cdot \operatorname{diag}(a_1^{-1},a_2^{-1},1)$ with $(\hat{x},\hat{y},\hat{z})$ substituted using $(\hat{x},\hat{y},\hat{z}) = (x,y,z) \operatorname{diag}(a_1,a_2,1) (S^T)^{-1}$.
\end{proof}

\begin{remark}
    The matrix $M$ from the proof of Theorem~\ref{thm:hedge2} appears in a natural way if one tries to find it as
    \[
        M=
        \begin{pmatrix}
            \hat{x} & \hat{y} & \hat{z} \\
            1+A_1 \hat{z} & B_1\hat{z} & C_1\hat{x} + D_1\hat{y}\\
            A_2 \hat{z} & 1+B_2\hat{z} & C_2 \hat{x} +D_2 \hat{y}
        \end{pmatrix}.
    \]
\end{remark}

\begin{remark}
    Let
    \begin{gather*}
    A:=k[a^{\pm 1}_1,a^{\pm 1}_2,a^{\pm 1}_3,x,y,z,x_0,x_1,x_2,x_3]/(q-1, \phi+a_1a_2a_3),\\\
    B:=k[a^{\pm 1}_1,a^{\pm 1}_2,a^{\pm 1}_3,x_0,x_1,x_2,x_3]/( \phi+a_1a_2a_3),
    \end{gather*}
    where
    \[
    q:=a_1x^2+a_2y^2+a_3z^2,\quad \phi:=x_0^2+a_2a_3x_1^2+a_1a_3x_2^2+a_1a_2a_3^2x_3^2.
    \]
    Define $\T$ to be the projective module over $A$ given as the kernel
    \[
    0\to \T \to A^{ 3} \xrightarrow{\nu} A \to 0,\quad \nu(u_1,u_2,u_3):=a_1xu_1+a_2yu_2+a_3zu_3.
    \]
    The inclusion of algebras $B\subseteq A$ corresponds to a fibration $\pi\colon\Spec A \to \Spec B$, and $\T$ is the relative tangent bundle for the fibration. The fibers of $\pi$ are affine quadrics $\{q=1\}$ together with a choice of a solution for the equation $\{\phi=-a_1a_2a_3\}$, so Theorem~\ref{thm:hedge2} yields that for $x\in \Spec B$ the vector bundle $\T|_{\pi^{-1}(x)}$ is trivial over the quadric $\pi^{-1}(x)$ (note that for this implication of the theorem there is no need to assume that the base field is perfect). It would be interesting to investigate the global structure of $\T$. 
\end{remark}

\section{Recollection on $\K$-cohomology and its variants}
For the investigation of nowhere vanishing sections and triviality of vector bundles we use cohomological obstructions as in \cite[Chapter~8]{Mor12} and \cite{Fas11b, Fas15} valued in $\H^\star(X,\K^\MW_*)$ and related cohomology groups. In this section we recall some of the relevant definitions and basic properties, and a detailed exposition may be found in \cite{Mor12,Fas19,Fel20,Fel21,DFJ25,Deg23}. 

\begin{recollection} 
Recall that the \textit{Milnor--Witt $\K$-theory} $\K^\MW_*(k)$ may be defined by generators $[a]$, $a\in k^\times$, of degree $1$, and $\eta$ of degree $-1$, and certain relations. For $a\in k^\times$ one puts $\langle a\rangle:=1+\eta [a] \in \K^\MW_0(k)$, and this induces an isomorphism $\K^\MW_0(k)\cong \GW(k)$ with $\langle a\rangle \in \K^\MW_0(k)$ corresponding to the quadratic form $\langle a\rangle$. Sending $\eta\mapsto 0$, $[a]\mapsto \{a\}$ gives a surjective homomorphism $\K^\MW_*(k)\to \K^\M_*(k)$ to Milnor $\K$-theory, while sending $h:=1+\langle-1\rangle\mapsto 0$, $[a]\mapsto \langle\langle a\rangle\rangle$, gives a surjective homomorphism $\K^\MW_*(k)\to\I^*(k)$, where $\I^{n}(k):=\I(k)^n$, $n\ge 0$, the respective power of the \textit{fundamental ideal} $\I(k):=\ker(\W(k)\xrightarrow{\dim} \Z/2)$, and $\I^n(k)=\W(k)$, $n\le 0$, is the Witt ring of regular quadratic forms over $k$. We also use the notation $\kM_*(k):=\K^\M_*(k)/2$ for the mod $2$ Milnor $\K$-theory and $\overline{\I}{}^*(k):=\I^*(k)/{\I}^{*+1}(k)$, and note that $\kM_*(k)\cong \overline{\I}{}^*(k)$ by Milnor conjecture \cite[Theorem~4.1]{OVV07}.

The graded rings $\K^\MW_*(k)$, $\K^\M_*(k)$, $\kM_*(k)$ and $\I^*(k)$ may be canonically extended to unramified sheaves on $\Sm_k$ keeping the same values on the function fields $k(X)$, $X\in \Sm_k$, allowing to consider $\H^n(X,M_*):=\H^n_{\Zar}(X,M_*)$ for $M_*=\K^\MW_*,\K^\M_*,\kM_*,\I^*$. Furthermore, for a line bundle $\mc{L}$ over $X\in \Sm_k$ one may introduce twisted cohomology groups $\H^n(X,M_*\otimes \mc{L})$, which admit canonical isomorphisms
\[
\H^n(X,M_*\otimes \mc{L}\otimes (\mc{L}')^2) \cong \H^n(X,M_*\otimes \mc{L}), \quad \mc{L}'\in\operatorname{Pic} X,
\]
and for $M_*$ being $\K^\M_*$ or $\kM_*$ canonical isomorphisms 
\[
\H^n(X,M_*\otimes \mc{L})\cong \H^n(X,M_*).
\]
We denote by 
\begin{gather*}
\widetilde{\CH}{}^n(X,\mc{L}):=\H^n(X,\KMW_n\otimes \mc{L}), \, \widetilde{\CH}{}^n(X):=\widetilde{\CH}{}^n(X,\struct_X),\, \widetilde{\CH}_n(X):=\widetilde{\CH}{}^{\dim X-n}(X,\omega_{X/k}),\\
\CH^n(X):=\H^n(X,\KM_n),\quad \Ch^n(X):=\H^n(X,\kM_n)\cong \CH^n(X)/2,
\end{gather*}
the (twisted) \textit{Chow--Witt groups}, \textit{Chow groups} and \textit{Chow groups modulo $2$}. 

The unramified graded sheaves $\K^\MW_*$, $\K^\M_*$, $\kM_*$ and $\I^*$ are particular examples of \textit{homotopy modules}, which admit an alternative description as \textit{Milnor--Witt cycle modules}, giving rise to a variant of the Gersten complex called the \textit{Rost--Schmid complexes}. The homology groups $A_n(X,M_{n+*}\otimes \mc{L})$ of the Rost--Schmid complex are usually referred to as (twisted) $n$-cycles with coefficients in~$M_*$ and enjoy an isomorphism
\[
A_n(X,M_{n+*}\otimes \mc{L}) \cong H^{d-n}(X,M_{d-n+*}\otimes (\mc{L}\otimes \omega_{X/k})),\quad d:=\dim X.
\]
We are mostly interested in top cohomology groups, equivalently, classes of (twisted) $0$-cycles, which admit a presentation
\begin{multline*}
H^{d}(X,M_{d+*}\otimes (\mc{L} \otimes  \omega_{X/k})) \cong A_0(X,M_{*}\otimes \mc{L}) = \\
=\coker \left(\bigoplus_{y\in X_{(1)}} M_{*+1}(k(y))\otimes (\mc{L}_y\otimes \omega_{y/k}) \xrightarrow{\partial} \bigoplus_{x\in X_{(0)}} M_*(k(x))\otimes \mc{L}_x\right),
\end{multline*}
where $X_{(1)}$ and $X_{(0)}$ are the sets of dimension $1$ and closed points of $X$ respectively, $\mc{L}_y$ and $\mc{L}_x$ are the fibers of $\mc{L}$ over $y$ and $x$ respectively, and $\omega_{y/k}$ is the canonical line bundle of $y$ over $\Spec k$. For a finite subset $Z\subseteq X_{(0)}$, an element $\alpha\in M_*(Z)$ and $\theta \in \mc{L}_Z$ we denote by $\alpha_Z\otimes \theta$ the corresponding twisted $0$-cycle supported on $Z$. Abusing the notation, we refer as a \textit{cycle} both to $\alpha_Z\otimes \theta$ and its class in the cohomology group.

For a finite separable field extension $L/k$ and a homotopy module $M_*$ we have a pushforward homomorphism $(\pi_{L/k})_*\colon M_*(L)\to M_*(k)$ which is functorial in $M_*$. For $M_*=\K^\M_*$ the pushforward homomorphism is induced by the norm homomorphism $\{a\}\mapsto \{N_{L/k}(a)\}$, while for $M_*=\I^*$ it is induced by Scharlau transfer. For $M_*=\K^\MW_*$ one can obtain the formula for the transfer combining the norm homomorphism with Scharlau transfer using the inclusion $\K^\MW_*\subseteq \K^\M_*\oplus \I^*$. In particular, for a quadratic field extension $k(\sqrt{\alpha})/k$ we have
\begin{equation}\label{eq:K1MW_transfer}
p_*([u+\sqrt{\alpha}v]) = \begin{cases}
    \langle -2u\rangle[-\alpha(u^2-\alpha v^2)] - \langle -2\rangle [-\alpha], &u\neq 0,\\
    \langle 2\rangle[-\alpha v^2], & u=0.
\end{cases}
\end{equation}

Cohomology groups with coefficients in a (twisted) homotopy module enjoy the usual properties of a cohomology theory including pullbacks $f^*\colon \H^n(X,M_*\otimes\mc{L}) \to \H^n(Y,M_*\otimes f^*\mc{L})$ for arbitrary morphisms $f\colon Y\to X$, pushforwards $f_*\colon A_n(Y,M_*\otimes f^*\mc{L}) \to A_n(X,M_*\otimes \mc{L})$ for proper morphisms, projection formula, localization long exact sequences and homotopy invariance. Furthermore, the group of classes of $0$-cycles with coefficients in a homotopy module is a birational invariant of a smooth proper variety by \cite[Theorem~2.2.12]{Fel22}, and in view of the following lemma, a stable birational invariant. The lemma below is certainly well known, but we could not find this precise statement in the literature and provide the proof for completeness.
\end{recollection}

\begin{lemma} \label{lem:push_proj_iso}
    For $X\in \Sm_k$ and a homotopy module $M_*$ the pushforward 
    \[
    \pi_*\colon A_0(\PP^n\times X,M_*) \to A_0(X,M_*)
    \]
    for the projection $\pi\colon \PP^n\times X\to X$ is an isomorphism.
\end{lemma}
\begin{proof}
    Consider a linear embedding $\PP^{n-1}\to \PP^n$ and the associated localization sequence
    \[
    A_0(\PP^{n-1}\times X,M_*) \xrightarrow{(i_n)_*} A_0(\PP^{n} \times X,M_*) \to A_0(\A^{n} \times X,M_*) \to 0
    \]
    Homotopy invariance yields $A_0(\A^{n} \times X,M_*)\cong A_{-n}(X,M_{*-n})=0$, hence $(i_n)_*$ is surjective. Applying induction we obtain that for the embedding $i\colon X \to \PP^n\times X$ given by the choice of a rational point the homomorphism
    \[
    i_*\colon A_0(X,M_*) \xrightarrow{} A_0(\PP^{n} \times X,M_*)
    \]
    is surjective. Since $i$ is a section for the projection $\pi$ then $i_*$ is also injective, thus $i_*$ is an isomorphism and $\pi_*$ is an isomorphism as well.
\end{proof}

\begin{corollary} \label{cor:stable_birational_invariant}
    Let $X\in \Sm_k$ be proper and rational,  and $M_*$ be a homotopy module. Then for the structure morphism $\pi\colon X\to \Spec k$ the pushforward homomorphism
    \[
    \pi_*\colon  A_0(X,M_*) \to M_*(k)
    \]
    is an isomorphism.
\end{corollary}
\begin{proof}
    Follows from \cite[Theorem~2.2.12]{Fel22} and Lemma~\ref{lem:push_proj_iso}.
\end{proof}

\begin{lemma} \label{lem:image_transfer_KMW1}
    Let $q$ be a regular quadratic form over $k$, $Q$ be the associated projective quadric, $n:=\dim Q\ge 0$, and
    \[
    \pi_*\colon \H^n(Q,\K^\M_{n+1})\to \K^\M_1(k)=k^\times,\quad \tilde{\pi}_*\colon \H^n(Q,\K^\MW_{n+1}\otimes \omega_{Q/k})\to \K^\MW_1(k),
    \]
    be the pushforwards along the structure morphism. Then for $a\in k^\times$ the following are equivalent:
    \begin{enumerate}
        \item $a\in [D(q)^2]$,
        \item $a\in\pi_*(\H^n(Q,\K^\M_{n+1}))$,
        \item $[a]\in\tilde{\pi}_*(\H^n(Q,\K^\MW_{n+1}\otimes \omega_{Q/k}))$.
    \end{enumerate}
\end{lemma}
\begin{proof}
    The equivalence (1)$\Leftrightarrow$(2) is \cite[Lemma~2.2]{CTS93} (see also \cite[Lemma~8]{Ros88}). The implication (3)$\Rightarrow$(2) follows from the commutativity of the square
    \[
    \xymatrix{
    \H^n(Q,\K^\MW_{n+1}\otimes \omega_{Q/k}) \ar[r]\ar[d]_{\tilde{\pi}_*} & \H^n(Q,\K^\M_{n+1})\ar[d]^{\pi_*} \\
    \K^\MW_1(k) \ar[r] & k^\times
    }
    \]
    where the horizontal maps are induced by the quotient morphism $\K^\MW_*\to \K^\M_*$. 
    
    (1)$\Rightarrow$(3). If $Q(k)\neq \emptyset$ then $\tilde{\pi}_*$ is surjective and there is nothing to do, so we may assume that $q$ is anisotropic. Rescaling $q$ we may assume that $1\in D(q)$, so we need to show that for $a\in [D(q)]$ we have $[a]\in\image \tilde{\pi}_*$. Pick $a_1,a_2\in k^\times$, by \cite[Lemma~3.5.1]{Mor12} we have $[a_1]+\langle a_1\rangle [a_2] =[a_1a_2]$. Since $\image \tilde{\pi}_*$ is stable under multiplication with the elements of $\K^\MW_{0}(k)$ by the projection formula, we have $[a_1a_2]\in \image \tilde{\pi}_*$ provided that $[a_1],[a_2]\in \image \tilde{\pi}_*$. Hence it suffices to show that $[a]\in \image \tilde{\pi}_*$ for $a\in D(q)$. 
    
    For the rest of the argument we follow the proof of \cite[Lemma~3.10]{AL25}.
    
    First suppose that $n=0$, then one may assume $q=x^2-\alpha y^2$ and $Q\cong \Spec k[\sqrt{\alpha}]$. We need to show that $[N_{k[\sqrt{\alpha}]/k}(u)]\in \image \tilde{\pi}_*$ for every $u\in k[\sqrt{\alpha}]^\times$. This easily follows from formulas~\eqref{eq:K1MW_transfer} and the fact that $\image \tilde{\pi}_*$ is a $\KMW_0(k)$-submodule of $\KMW_1(k)$.

    Now suppose $n\ge 1$. The same argument as in the end of the proof of \cite[Lemma~3.10]{AL25} shows that there is $c\in D(q)$ and closed points $z_1,z_2\in Q$ such that $k(z_1)\cong k[\sqrt{-ac^{-1}}]$, $k(z_2)\cong k[\sqrt{-c^{-1}}]$. Then
    \[
    \tilde{\pi}_*([\sqrt{-ac^{-1}}]_{z_1})=\langle2\rangle[ac^{-1}],\quad \tilde{\pi}_*([\sqrt{-c^{-1}}]_{z_2})= \langle2\rangle[c^{-1}].
    \]
    Applying \cite[Lemma~3.5]{Mor12} we get
    \[
    [ac^{-1}]=[a]-\langle ac^{-1} \rangle [c],\quad [c^{-1}]=-\langle c^{-1}\rangle [c].
    \]
    Since $\image \tilde{\pi}_*$ is a $\K^\MW_{0}(k)$-submodule of $\K^\MW_{1}(k)$ and $[ac^{-1}],[c^{-1}]\in \image \tilde{\pi}_*$ by the above, we have $[a]\in \image \tilde{\pi}_*$ and the claim follows.
\end{proof} 

\section{A digression on Quillen $\K$-theory of quadrics and higher Chern classes}

In this section $q=\sum_{i=1}^{n+1}a_ix_i^2$ is a regular quadratic form over $k$ and we put
\[
\begin{gathered}
    \tilde{q}:=-x_0^2+q,\quad Q:=\{\tilde{q}=0\}\subseteq \PP^{n+1},\quad Q^\infty:=Q\cap \{x_{0}=0\},\quad Q^o:=Q\setminus Q^\infty.
    \end{gathered}
\]
Below we study the properties of the Quillen $\K$-theory groups $\K_1(Q)$ and $\K_1(Q^o)$ and higher Chern class 
    \[
        c_{1,n+1}:=\Omega c_{n+1}\colon \K_1(Q)\to \H^{n}(Q,\K^\M_{n+1})
    \]
defined as in \cite{Pus04}. These will be used in the further sections to investigate cohomological obstructions to the triviality of the tangent bundle $T_{Q^o}$. For a ring $A$ and $u\in A^\times$ we denote by $\{u\}$ the associated element of $\K_1(A)$.

\begin{lemma} \label{lem:K-theory_quadrics_surjection}
    Suppose that either $n$ is even or $n=2m+1$ and $\Delta:=(-1)^{m+1} \prod_{i=1}^{2m+2} a_i \notin (k^\times)^2$. Then the pullback homomorphism $\K_1(Q)\to \K_1(Q^o)$ is surjective.
\end{lemma}
\begin{proof}
     Consider the following part of the localization sequence for Quillen $\K$-theory.
    \[
    \K_1(Q) \xrightarrow{} \K_1(Q^o)\xrightarrow{} \K_0(Q^\infty)\xrightarrow{i_*} \K_0(Q)
    \]
    The claim is equivalent to injectivity of $i_*$. Since \cite[Theorem~9.1]{Swa85} yields that $\K_0(Q^\infty)$ is a free abelian group, it is sufficient to check injectivity rationally. Using the Chern character isomorphism and the Riemann-Roch theorem it suffices to prove that 
    \[
    i_*\colon \CH^*(Q^\infty)_\QQ \to \CH^*(Q)_\QQ
    \]
    is injective. Consider base change to the algebraic closure.
    \[
    \xymatrix{
    \CH^*(Q^\infty)_\QQ \ar[r]^{i_*}\ar[d]_{p_\infty^*} & \CH^*(Q)_\QQ \ar[d]^{p^*}\\
    \CH^*(Q_{\bar{k}}^\infty)_\QQ \ar[r]^{\bar{\imath}_*}  & \CH^*(Q_{\bar{k}})_\QQ
    }
    \]
    Here both vertical maps are injective, we refer to \cite[Proposition~1.1]{KRS98} for this and the properties used below.

    If $n=2m+1$ then
    \[
    \CH^*(Q_{\bar{k}}^\infty)_\QQ=\bigoplus_{j=0}^{m} \QQ\cdot [h]^j \oplus \bigoplus_{j=0}^{m}\QQ\cdot [\ell_j], \quad \CH^*(Q_{\bar{k}})_\QQ=\bigoplus_{j=0}^{m} \QQ\cdot [h]^j \oplus \bigoplus_{j=0}^{m}\QQ\cdot [\ell_j],
    \]
    where $[h]$ is the class of a hyperplane section, $[\ell_j]$ is the class of a dimension $j$ linear subspace of a quadric. Since $Q^\infty$ is a hyperplane section we have 
    \[
    \bar{\imath}_*([h]^j)=[h]^{j+1},\,j<m,\quad \bar{\imath}_*([h]^m)=[h]^{m+1}=2[\ell_m],\quad \bar{\imath}_*([\ell]_j)=[\ell_{j}].
    \]
    Since the signed discriminant is not a square by assumption,
    \[
    \image p_\infty^*=\bigoplus_{j=0}^{m} \QQ\cdot [h]^j \oplus \bigoplus_{j=0}^{m-1}\QQ\cdot [\ell_j].
    \]
    The above formulas yield injectivity of $\bar{\imath}_*p^*_\infty=p^*i_*$, hence injectivity of $i_*$ giving the claim of the lemma. The proof for $n$ even is similar but easier, since in this case $\bar{\imath}_*$ is injective.
\end{proof}

\begin{definition}
    The $n$-th \textit{power sum Chern class} $s_n(\mc{V})\in \CH^n(X)$ for a vector bundle $\mc{V}$ over $X$ is the characteristic class satisfying (i) $s_n(\mc{L})=c_1(\mc{L})^n$ for a line bundle $\mc{L}$ over $X$ and (ii) $s_n(\mc{V}_1\oplus \mc{V}_2) = s_n(\mc{V}_1) + s_n(\mc{V}_2)$ for vector bundles $\mc{V}_1, \mc{V}_2$ over $X$. It is uniquely defined by the splitting principle and can be computed from Chern classes using the polynomial expressing the power sum $\sum x_i^n$ via elementary symmetric polynomials $\sigma_1,\hdots, \sigma_n$.
\end{definition}

\begin{lemma} \label{lem:higher_Chern_formula}
    Let $X\in \Sm_k$, $\mc{V}$ be a vector bundle over $X$ and $u\in k^\times$. Put $n=\dim X$. Then
    \[
    c_{1,n+1} (\{u\}\cdot [\mc{V}])= (-1)^n\{u\}\cdot s_n(\mc{V}).
    \]
\end{lemma}
\begin{proof}
    Splitting principle together with the additivity of $c_{1,n+1}$ and $s_n$ yield that it suffices to prove the claim for $\mc{V}=\mc{L}$ a line bundle.
    By \cite[Proposition~4.1]{Shi14} we have
    \[
    c_{1,n+1}(\{u\}\cdot [\mc{L}]) = \sum_{i=0}^{n} (-1)^i {n \choose i} c_{1,n+1-i}(\{u\}) c_1(\mc{L})^i .
    \]
    By naturality and vanishing $\H^{n-i}(\Spec k, \K^\M_{n+1-i})=0$ one has $c_{1,n+1-i}(\{u\})=0$ if $n-i\neq 0$, and \cite[Remark~V.11.11.1]{Wei13} gives $c_{1,1}(\{u\})=\{u\}$, thus 
    \[
    c_{1,n+1}(\{u\}\cdot [\mc{L}]) =(-1)^n c_{1,1}(\{u\})\cdot c_1(\mc{L})^n=(-1)^n\{u\}\cdot c_1(\mc{L})^n. \qedhere
    \]
\end{proof}

\begin{lemma} \label{lem:Nrd_in_d2}
    Let $n=3$. Then $\Nrd_{C_0(\tilde{q})/k}(C_0(\tilde{q})^\times)\subseteq [D(q)^2]$.
\end{lemma}
\begin{proof}
    It follows from \cite[Exercise~85.4]{EKM08} applied to $\varphi=\tilde{q}$ that if for a field extension $L/k$ the central simple algebra $C_0(\tilde{q})_L$ splits then $\tilde{q}_L$ is split, and hence $q_L$ is isotropic, being a codimension $1$ subform. Since by \cite[Proposition~2.6.8]{GS17} we have
    \[
    \Nrd_{C_0(\tilde{q})/k}(C_0(\tilde{q})^\times) = [\{N_{L/k}(a)\mid \text{$a\in L^\times$, $L\subseteq C_0(\tilde{q})$ and $C_0(\tilde{q})_L$ splits}\}],
    \]
    the claim follows from Lemma~\ref{lem:CTS}.
\end{proof}

\begin{lemma} \label{lem:K-theory_quadrics_push_higher}
    Consider the pushforward homomorphism $\pi_*\colon \H^n(Q,\K^\M_{n+1}) \to k^\times$. Suppose $n\ge 3$ is odd. Then the following holds.
    \begin{enumerate}
        \item If $n=3$ then $\pi_* c_{1,n+1}(\K_1(Q)) \subseteq [D(q)^2]$.
        \item If $n\ge 5$ then $\pi_* c_{1,n+1}(\K_1(Q))= (k^\times)^2\subseteq [D(q)^2]$.
    \end{enumerate}
\end{lemma}
\begin{proof}
    Let $A:=C_0(\tilde{q})$ be the even part of the Clifford algebra. Since $n$ is odd, $A$ is a central simple algebra over $k$ of dimension $N^2$, $N:=2^{\frac{n+1}{2}}$. Consider the isomorphism
    \[
    \K_1(Q)=\bigoplus_{i=0}^{n-1} \K_1(k)\cdot[\struct(-i)] \oplus \K_1(A)
    \]
    constructed in \cite[\S~10]{Pan94} (note that it slightly differs from the one constructed in \cite[Theorem~9.1]{Swa85}).
    Here the inclusion $\Phi\colon\K_1(A)\to \K_1(Q)$ is induced by the functor 
    \[
    A\lmod\to \Coh Q,\quad P\mapsto \mc{U}\otimes _{A} P,
    \]
    where $\mc{U}$ is a particular vector bundle over $Q$ equipped with a structure of a right $A$-module. It is straightforward to check that over the algebraic closure one has $\mc{U}_{\bar{k}}\cong \mc{S}\otimes_{\bar{k}} \bar{k}^{\oplus N}$ for the spinor bundle $\mc{S}$ over $Q_{\bar{k}}$ defined as in \cite[Definition~1.3]{Ott88}, and that the action of $A_{\bar{k}}\cong \M_N(\bar{k})$ on $\mc{U}_{\bar{k}}\cong \mc{S}\otimes_{\bar{k}} \bar{k}^{\oplus N}$ is given by the right action on $\bar{k}^{\oplus N}$.
    
    First we compute $\pi_* c_{1,n+1}(\{u\}\cdot [\struct(-i)])$ for $u\in k^\times$. By Lemma~\ref{lem:higher_Chern_formula} we have
    \[
    c_{1,n+1}(\{u\}\cdot [\struct(-i)]) =(-1)^n \{u\}\cdot c_1(\struct(-i))^n=i^n \{u\}\cdot c_1(\struct(1))^n.
    \]
    We have $\pi_*(c_1(\struct(1))^n)=2$, thus the projection formula yields
    \[
    \pi_*c_{1,n+1}(\{u\}\cdot [\struct(-i)]) = u^{2i^n}
    \]
    which is a square.

    It remains to compute $\pi_* c_{1,n+1}(\Phi(\alpha))$ for $\alpha\in \K_1(A)$ . In order to do so we consider the following commutative diagram for the base change to the algebraic closure:
    \[
    \xymatrix{
       \K_1(A)\ar[r]^{\Phi} \ar[d]_f & \K_1(Q) \ar[r]^(0.4){c_{1,n+1}} \ar[d] & \H^n(Q,\K^\M_{n+1})\ar[r]^(0.6){\pi_*} \ar[d] & k^\times \ar[d]^\iota  \\
       \K_1(A_{\bar{k}}) \ar[r]^{\Phi_{\bar{k}}} & \K_1(Q_{\bar{k}}) \ar[r]^(0.4){c_{1,n+1}} & \H^n(Q_{\bar{k}},\K^\M_{n+1})\ar[r]^(0.6){\pi_*} & \bar{k}^\times 
    }
    \]
    Over the algebraic closure the Clifford algebra splits as $A_{\bar{k}}\cong \M_{N}(\bar{k})$. Morita equivalence $\M_{N}(\bar{k})\lmod \simeq \bar{k}\lmod$ yields a commutative triangle
    \[
    \xymatrix{
       \K_1(A_{\bar{k}}) \ar[rr]^{\Phi_{\bar{k}}} \ar[rd]_{\Theta}^{\cong} & & \K_1(Q_{\bar{k}})  \\
       & \K_1(\bar{k}) \ar[ru]_{\Psi} &
    }
    \]
    where $\Theta (\bar{u})=\{\det \bar{u}\}$ for $\bar{u}\in A_{\bar{k}}^\times \cong \GL_N(\bar{k})$ and $\Psi(a)=\{a\}\cdot [\mc{S}_{\bar{k}}]$ for $a\in \bar{k}^\times\cong \K_1(\bar{k}) \subseteq \K_1(Q_{\bar{k}})$, and $[\mc{S}_{\bar{k}}]\in \K_0(Q_{\bar{k}})$. Then for $u\in A^\times$ we have
    \[
    \Theta f(\{u\}) = \{\det \bar{u}\} = \{\Nrd_{A/k} u\}
    \]
    with $\Nrd_{A/k} u\in k^\times \subseteq \bar{k}^\times$ being the reduced norm of $u$. Applying the above and Lemma~\ref{lem:higher_Chern_formula} we obtain 
    \[
    \iota\pi_*c_{1,n+1}\Phi(\{u\}) = \pi_*c_{1,n+1}(\{\Nrd_{A/k} u\}\cdot [\mc{S}_{\bar{k}}]) = (\Nrd_{A/k} u)^{-\deg s_n(\mc{S})}, \quad u\in A^\times.
    \] 
    Since $A^\times \to \K_1(A)$, $u\mapsto \{u\}$, is surjective, and $\iota$ is injective, the claim for $n=3$ follows from Lemma~\ref{lem:Nrd_in_d2}. It remains to show that $\deg s_n(\mc{S})$ is even for $n\ge 5$. Recall \cite[Definition~1.3 and Theorem~2.8.(i)]{Ott88} that $\rank \mc{S}=\frac{N}{2}$, $\mc{S}^\vee\cong \mc{S}(1)$ and there is a short exact sequence
    \[
    0\to \mc{S}\to \struct_{Q_{\bar{k}}}^{\oplus N} \to \mc{S}(1)\to 0.
    \]
    Put $h:=c_1(\struct_{Q_{\bar{k}}}(1))$ and $d_i:=\deg h^{n-i}s_i$. Then from $\deg h^n=2$ we have $d_0=N=2^{\frac{n+1}{2}}$. The short exact sequence gives an equality $[\mc{S}]+[\mc{S}(1)]=N$ in $\K_0(Q_{\bar{k}})$, which combined with $[\mc{S}^\vee] = [\mc{S}(1)]$ yields that $s_i(\mc{S})=0$ for $i>0$ even, hence $d_i=0$ for $i>0$ even. For $i$ odd from $[\mc{S}^\vee] = [\mc{S}(1)]$ we obtain a recursive formula
    \[
    2d_{2m-1} = -2^{\frac{n+1}{2}} - \sum_{i=1}^{m-1} {2m-1 \choose 2i-1} d_{2i-1}.
    \]
    Hence we have $d_1=-2^{\frac{n-1}{2}}$, $d_3=2^{\frac{n-3}{2}}$, $d_5=-2^{\frac{n-1}{2}}$. Induction and recursive formula yield that $v_2(d_{2m-1})\ge \frac{n+1}{2}+1-m$ for $3\le m\le \frac{n+1}{2}$, where $v_2$ is the $2$-adic valuation. Hence $v_2(d_{n})\ge 1$ giving the claim.
\end{proof}

\section{Hermitian $\K$-theory of quadrics}

Using Swan's decomposition of the category of perfect complexes over a smooth quadric $Q$, H.~Xie was able to compute the Hermitian $\K$-theory of $Q$ in \cite{Xie26}. We gather here a few salient points that will prove useful in the rest of the paper. First, recall that the (twisted) Hermitian $\K$-theory spaces of $Q$ are spectra of the form $\mathrm{GW}^{[d]}(Q,\mathcal L)$, where $[d]\in\Z/4\Z$ and $\mathcal L$ is a line bundle over $Q$. There are canonical equivalences
\[
\mathrm{GW}^{[d]}(Q,\mathcal L\otimes \mathcal{N}^{\otimes 2})\simeq \mathrm{GW}^{[d]}(Q,\mathcal L)
\]   
for any line bundle $\mathcal{N}$ over $Q$. The Clifford algebra $C(q)$ is equipped with the so-called canonical involution $\sigma\colon C(q)\xrightarrow{\simeq} C(q)^{\mr{op}}$, which induces involutions on both the even part $C_0(q)$ and the odd part $C_1(q)$ of the algebra. Note that both $C_0(q)$ and $C_1(q)$ are $C_0(q)$-bimodules, and that the involution allows to consider for any left $C_0(q)$-module $M$ the right $C_0(q)$-module $\mathrm{Hom}_{C_0(q)}(M,C_i(q))$ as a left $C_0(q)$-module via the formula $(\lambda\cdot f)(m)=f(m)\sigma(\lambda)$ for any $\lambda \in C_0(q)$, $m\in M$ and $f\in \mathrm{Hom}_{C_0(q)}(M,C_i(q))$. 

Suppose that $q=\langle a_1,\ldots,a_{n+2}\rangle$ for some $n\in\mathbb{N}$, and write $e_1,\ldots,e_{n+2}$ for the canonical basis of $k^{n+2}$. If $n$ is odd, then $e:=e_1\cdot\ldots\cdot e_{n+2}\in C_1(q)$ is central in $C(q)$ and satisfies $\sigma(e)=(-1)^{\frac {n+1}2}e$. As $e\sigma(e)=\prod_{i=1}^{n+2}a_i\in k^\times$, it follows that 
\[
r_e\colon C_0(q)\to C_1(q), \phantom{i}x\mapsto xe
\]
induces an isomorphism of $C_0(q)$-bimodules such that the following diagram
\[
\xymatrix@C=5em{C_0(q)\ar[r]^-{\sigma}\ar[d]_-{r_e} & C_0(q)\ar[d]^-{r_e} \\
C_1(q)\ar[r]_-{(-1)^{\frac {n+1}2}\sigma} & C_1(q)}
\]
commutes. For any left $C_0(q)$-module $M$, we obtain an isomorphism of left $C_0(q)$-modules 
\[
\psi_M\colon M^{\sharp_0}:=\mathrm{Hom}_{C_0(q)}(M,C_0(q))\to \mathrm{Hom}_{C_0(q)}(M,C_1(q))=:M^{\sharp_1}, \phantom{i}f\mapsto r_e\circ f.
\]
Consider next the morphism of left $C_0(q)$-modules 
\[
\omega_M^i\colon M\to M^{\sharp_i\sharp_i}
\]
defined by $m\mapsto m^*$ with $m^*(f)=\sigma(f(m))$ for any $f\in M^{\sharp_i}$. Using the above diagram, it is straightforward to check that we have a commutative diagram
\[
\xymatrix{M\ar[r]^-{\omega_M^0}\ar[d]_-{(-1)^{\frac {n+1}2}\omega_M^1} & M^{\sharp_0\sharp_0}\ar[d]^-{\psi_{M^{\sharp_0}}} \\
M^{\sharp_1\sharp_1}\ar[r]_{\psi_M^{\sharp_1}} & (M^{\sharp_0})^{\sharp_1}}
\]
inducing an equivalence of spectra for $n=1\pmod 4$
\[
\mathrm{GW}^{[d]}(C_0(q)_{\sigma},C_1(q)_{\sigma})\simeq \mathrm{GW}^{[d+2]}(C_0(q)_{\sigma})
\]
and an equivalence 
\[
\mathrm{GW}^{[d]}(C_0(q)_{\sigma},C_1(q)_{\sigma})\simeq \mathrm{GW}^{[d]}(C_0(q)_{\sigma})
\]
for $n=3\pmod 4$. Here, we follow the notation of \cite{BP05} and \cite{Xie19}, stressing that the objects considered depend on the involution $\sigma$ and that $C_1(q)_{\sigma}$ is a duality coefficient in the sense of \cite[\S 7.4]{Sch10b}. 
We are now ready to state the main result of \cite{Xie26} under the form needed below.

\begin{theorem}\label{thm:GWspacesquadrics}
If $Q$ is a smooth projective quadric of dimension $n\geq 1$, there is for any $[d]\in\Z/4\Z$ a fibre sequence of spectra of the form
\[
\mathrm{K}(k)^{\lceil \frac n2-1\rceil}\oplus \mathrm{GW}^{[d]}(\mathfrak A(\omega_{Q/k}))\to \mathrm{GW}^{[d]}(Q,\omega_{Q/k})\xrightarrow{p_*}\mathrm{GW}^{[d-n]}(k)
\]
where $p_*$ is the pushforward map induced by $p\colon Q\to\Spec(k)$ and the space $\mathrm{GW}^{[d]}(\mathfrak A(\omega_{Q/k}))$ admits the following description:
\begin{enumerate}
\item $\mathrm{GW}^{[d]}(\mathfrak A(\omega_{Q/k}))\simeq \mathrm{GW}^{[d+2]}(C_0(q)_{\sigma})$ if $n=1\pmod 4$.
\item $\mathrm{GW}^{[d]}(\mathfrak A(\omega_{Q/k}))\simeq \mathrm{GW}^{[d]}(C_0(q)_{\sigma})$ if $n=3\pmod 4$.
\item If $n$ is even, there is a triangle of the form
\[
\mathrm{GW}^{[d]}(k)\to \mathrm{GW}^{[d]}(\mathfrak A(\omega_{Q/k}))\to \mathrm{GW}^{[d+1]}(C_0(q)_{\sigma})\xrightarrow{\eta\cup\mathrm{tr}}\mathrm{GW}^{[d]}(k)[1]
\]
where $\eta$ is the Bott element and $\mathrm{tr}$ is a convenient trace map.
\end{enumerate} 
\end{theorem}

\begin{proof}
Combine \cite[Theorem 1.1]{Xie26} with the isomorphism $\omega_{Q/k}\simeq \mathcal{O}(-n)$ and the above computation of $\mathrm{GW}^{[d]}(C_0(q)_{\sigma},C_1(q)_{\sigma})$ for odd $n$.
\end{proof}

\begin{corollary}\label{cor:exactsequences}
If $n$ is odd or $n=2\pmod 4$, the pushforward homomorphism
\[
p_*\colon \mathrm{W}^n(Q,\omega_{Q/k})\to \mathrm{W}(k)
\]
is injective. If $n=0\pmod 4$, there is an exact sequence
\[
\mathrm{W}^0(C_0(q)_{\sigma})\xrightarrow{\mathrm{tr}} \mathrm{W}(k)\xrightarrow{p^*}\mathrm{W}^0(Q,\omega_{Q/k})\xrightarrow{p_*}\mathrm{W}(k)
\]
\end{corollary}

\begin{proof}
If $n$ is odd, we apply $\pi_{-4}$ to the fiber sequence
\[
\mathrm{K}(k)^{\lceil \frac n2-1\rceil}\oplus\mathrm{GW}^{[n]}(\mathfrak A(\omega_{Q/k}))\to \mathrm{GW}^{[n]}(Q,\omega_{Q/k})\xrightarrow{p_*}\mathrm{GW}^{[0]}(k)
\]
to obtain an exact sequence
\[
\pi_{-4}\mathrm{GW}^{[n]}(\mathfrak A(\omega_{Q/k}))\to \pi_{-4}\mathrm{GW}^{[n]}(Q,\omega_{Q/k})\xrightarrow{p_*}\pi_{-4}\mathrm{GW}^{[0]}(k).
\]
Now $\pi_{-4}\mathrm{GW}^{[n]}(\mathfrak A(\omega_{Q/k}))=\mathrm{W}^{[n+4]}(\mathfrak A(\omega_{Q/k}))=\mathrm{W}^{[n]}(\mathfrak A(\omega_{Q/k}))$ and the same holds for the other groups. We are thus left to show that $\mathrm{W}^{[n]}(\mathfrak A(\omega_{Q/k}))=0$. If $n=1\pmod 4$, the equivalence
\[
\mathrm{GW}^{[n]}(\mathfrak A(\omega_{Q/k}))\simeq \mathrm{GW}^{[n+2]}(C_0(q)_{\sigma})
\]
yields an isomorphism $\mathrm{W}^{[n]}(\mathfrak A(\omega_{Q/k}))\simeq \mathrm{W}^{[n+2]}(C_0(q)_{\sigma})$, and the latter vanishes by \cite[Theorem A.3]{Xie19}. The same holds if $n=3\pmod 4$, using this time the equivalence 
\[
\mathrm{GW}^{[n]}(\mathfrak A(\omega_{Q/k}))\simeq \mathrm{GW}^{[n]}(C_0(q)_{\sigma}).
\]
If now $n$ is even, the result follows directly from \cite[Theorem 1.2]{Xie19} using $\mathrm{W}^{i}(k)=0$ for $i\neq 0\pmod 4$.
\end{proof}

This corollary allows to understand the kernel of the pushforward homomorphism 
\[
p_*^{\mathrm{W}}\colon \H^n(Q,\mathrm{W}\otimes \omega_{Q/k})\to \mathrm{W}(k)
\]
in low dimensions, as shown by the next result.

\begin{theorem}\label{thm:explicitkernel}
Suppose that $Q$ is a smooth projective quadric of dimension $1\leq n\leq 7$. If $n\neq 4$, there is an exact sequence
\[
\mathrm{W}^{n-5}(Q,\omega_{Q/k})\xrightarrow{e}\H^{n-5}(Q,\mathrm{W}\otimes\omega_{Q/k})\xrightarrow{\partial} \H^{n}(Q,\mathrm{W}\otimes\omega_{Q/k})\xrightarrow{p_*^{\mathrm{W}}}\mathrm{W}(k)
\]
where $\partial$ is the differential at page $5$ of the Gersten--Witt spectral sequence and $e$ is an edge map. In particular, $p_*^{\mathrm{W}}$ is injective if $1\leq n\leq3$. If $n=4$, then there is an exact sequence of the form
\[
0\to \ker(sp^*)/\ker(p^*)\to \H^{4}(Q,\mathrm{W}\otimes\omega_{Q/k})\xrightarrow{p_*^{\mathrm{W}}}\mathrm{W}(k),
\]
where $sp^*\colon \mathrm{W}(k)\to  \mathrm{W}(k(Q))$ is induced by $k\subset k(Q)$ and $p^*\colon \mathrm{W}(k)\to \mathrm{W}(Q,\omega_{Q/k})$ is the pullback.
\end{theorem}

\begin{proof}
Suppose first that $n\neq 4$. In that case, we know from Corollary \ref{cor:exactsequences} that the pushforward map $p_*\colon \mathrm{W}^n(Q,\omega_{Q/k})\to \mathrm{W}(k)$ is injective, and we are left to understand the kernel of the edge map $\H^n(Q,\mathrm{W}\otimes\omega_{Q/k})\to \mathrm{W}^n(Q,\omega_{Q/k})$ in the Gersten--Witt spectral sequence. This is a direct application of \cite[Theorem 10.1]{BW02}. 

If $n=4$, we have an exact sequence $\mathrm{W}(k)\xrightarrow{p^*}\mathrm{W}^0(Q,\omega_{Q/k})\xrightarrow{p_*}\mathrm{W}(k)$ and the Gersten--Witt spectral sequence yields an exact sequence
\[
0\to \H^4(Q,\mathrm{W}\otimes\omega_{Q/k})\xrightarrow{j} \mathrm{W}^0(Q,\omega_{Q/k})\to \H^0(Q,\mathrm{W}\otimes\omega_{Q/k})\to 0.
\]
We observe that $\mathrm{W}(k)\to  \mathrm{W}(k(Q))$ factors as $\mathrm{W}(k)\to \H^0(Q,\mathrm{W}\otimes\omega_{Q/k})\to  \mathrm{W}(k(Q))$ with injective second map and we consider the following diagram
\[
\xymatrix{ & \H^4(Q,\mathrm{W}\otimes\omega_{Q/k})\ar[d]^-j\ar[rd]^-{p_*^{\mathrm{W}}} & \\
\mathrm{W}(k)\ar[r]^-{p^*}\ar[rd]_-{sp^*} &  \mathrm{W}^0(Q)\ar[r]^-{p_*}\ar[d]^-s &\mathrm{W}(k). \\
 & \mathrm{W}(k(Q)) & }
\]
A straightforward diagram chase yields an isomorphism $\ker(p_*^{\mathrm{W}})\to  \ker(sp^*)/\ker(p^*)$, using that $j$ is injective. 
\end{proof}

Finally, let us prove the following result that will prove useful in computations.

\begin{corollary}\label{cor:K1symplecticinjective}
If $Q\subset \mathbb{P}^2$ is a smooth conic, the pushforward homomorphism
\[
p_*\colon \mathrm{GW}_1^{[2]}(Q,\omega_{Q/k})\to \mathrm{GW}_1^{[1]}(k)
\] 
is injective.
\end{corollary}

\begin{proof}
As above, it suffices to prove that 
\[
\mathrm{GW}_1^{[2]}(\mathfrak A(\omega_{Q/k}))=\mathrm{GW}_1^{[0]}(C_0(q)_{\sigma})=0
\]
to conclude. As $n=1$, $A:=C_0(q)$ is a (possibly split) quaternion algebra, and $\sigma$ coincides with the quaternionic conjugation $\gamma$. Thus, we are left to show that $\mathrm{GW}_1^{[0]}(A,\gamma)=0$. With this in mind, observe that $\mathrm{GW}_1^{[0]}(A,\gamma)$ is the abelianization of the stable isometry group of hermitian hyperbolic forms (e.g. \cite[Remark 4.13]{Sch10}) which vanishes by \cite[Theorems 5.9 and 6.1(3)]{Gil09}.
\end{proof}

\section{Triviality in higher dimension: necessary and sufficient conditions}

In this section we establish some necessary and sufficient conditions for the tangent bundle $T_{Q^o}$ to be trivial. This is done by analyzing cohomological obstructions from \cite{Fas11b, Fas15}. 

\begin{recollection}
    Let $X$ be a smooth variety. We denote by $\Vect_n(X)$ the set of isomorphism classes of rank $n$ vector bundles over $X$ and 
    \[
    \Vect^o_{n}(X):=\{V\in \Vect_n(X)\mid V\oplus \struct_X \cong \struct_X^{n+1}\}.
    \]
    If $X=\Spec R$ then there is a bijection \cite[Proposition~1.4.8]{Lam06}
    \[
    \Um_{n+1}(R)/\GL_{n+1}(R)\cong \Vect^o_{n}(X),
    \]
    where a unimodular row is mapped to the kernel of the surjection $R^{n+1}\to R$ given by the row. We may also restrict the action of $\GL_{n+1}(R)$ to a subgroup and consider $\Um_{n+1}(R)/\mr{E}_{n+1}(R)$ for the elementary subgroup $\mr{E}_{n+1}(R)\leqslant \GL_{n+1}(R)$ generated by elementary transvections.

    Let $X:=\Spec R$ be a smooth affine variety and put $n:=\dim X$. Then following \cite[\S~4]{Fas11b} one may define a map 
    \[
    \phi\colon \Um_{n+1}(R) \to A_0(X,\K^\MW_{1}\otimes \omega^\vee_{X/k}) \cong \H^n (X,\K^\MW_{n+1}).
    \]
    If the zero subscheme $Z:=Z(u_2,u_3,\hdots,u_{n+1})\subseteq X$ is reduced and finite, then by the computation in \cite[\S~3.3]{Fas11b} one has
    \[
    \phi(u_1,u_2,\hdots,u_{n+1})=\sum_{x\in Z} [\bar{u}_1]_x\otimes (du_2\wedge \hdots \wedge du_{n+1})_x^\vee,
    \]
    where $[\bar{u}_1]_x\in \K_1^\MW(k(x))$ is the symbol associated with the image $\bar{u}_1\in k(x)$  of $u_1$ in the residue field, and $(du_2\wedge \hdots \wedge du_{n+1})_x^\vee\in (\omega_{X/k}^\vee)_x$ is the element dual to $(du_2\wedge \hdots \wedge du_{n+1})_x$. One may show \cite[Theorem~4.9]{Fas11b} that for $n\ge 2$ this induces a bijection
    \begin{equation} \label{eq:obstruction_general_triviality}
    \Um_{n+1}(R)/\mr{E}_{n+1}(R) \xrightarrow{\simeq} \H^n (X,\K^\MW_{n+1}).
    \end{equation}
    More precisely, the left-hand term is an abelian group according to \cite[Theorem 4.1]{vdK89}, and the map is an isomorphism. Note that \cite[Theorem~4.9]{Fas11b} needs the assumption that $n\geq 3$. However, one can use the formalism of Postnikov towers combined with \cite[Theorem~1]{AHW17} to see that it holds for $n=2$ as well.
    Composing $\phi$ with the quotient morphism we obtain a map
    \[
    \bar{\phi} \colon  \Um_{n+1}(R) \to A_0(X,\K^\M_{1}) \cong \H^n (X,\K^\M_{n+1}).
    \]
    If the zero subscheme $Z:=Z(u_2,u_3,\hdots,u_{n+1})\subseteq X$ is reduced and finite then 
    \[
    \bar{\phi}(u_1,u_2,\hdots,u_{n+1})=\sum_{x\in Z} \{\bar{u}_1\}_x
    \]
    with $\{\bar{u}_1\}\in \K_1^\M(k(x))=k(x)^\times$ being the image of $u_1$ in the residue field. If $n$ is odd this induces an isomorphism (\cite[Proposition 2.8]{Fas15})
    \begin{equation} \label{eq:obstruction_odd_triviality}
    \Um_{n+1}(R)/(\GL_{n+1}(R)\cap \mr{E}_{n+2}(R)) \xrightarrow{\simeq} \H^n(X,\K^\M_{n+1}).
    \end{equation}
    Associating to a matrix its first row yields a map $\GL_{n+1}(R)\to \Um_{n+1}(R)$, and in turn a map
    \[
    r\colon \GL_{n+1}(R)/(\GL_{n+1}(R)\cap E_{n+2}(R))\to \Um_{n+1}(R)/(\GL_{n+1}(R)\cap E_{n+2}(R)).
    \] 
    The left-hand side is actually the group $\mathrm{K}_1(R)$ (by stabilization), and the map is a group homomorphism by \cite[Theorem 5.3 (ii)]{vdK89}. In view of this, one concludes \cite[Remark 2.12]{Fas15} that for $n$ odd there is an exact sequence
    \begin{equation}\label{eq:stably_trivial_as_cokernel}
     \mathrm{K}_1(R)\xrightarrow{r}\H^n (X,\K^\M_{n+1})\to \Um_{n+1}(R)/\GL_{n+1}(R)\cong \Vect^o_n(X).
     \end{equation}
     Finally, note that the homomorphism $r$ was shown in \cite[Proposition~6.3]{Du22} (take $\xi=0$ in the statement) to coincide with the (simplicial) loop of the Chern class $c_{n+1}$, which is the higher Chern class $c_{1,n+1}$.
\end{recollection}

Below $q=\sum_{i=1}^{n+1}a_ix_i^2$ is a regular quadratic form over $k$, $\tilde{q}:=-x_0^2+q$ and we put
\[
\begin{gathered}
    Q:=\{\tilde{q}=0\}\subseteq \PP^{n+1},\quad Q^\infty:=Q\cap \{x_{0}=0\},\quad Q^o:=Q\setminus Q^\infty\cong \Spec R,\\
    R=k[x_1,x_2,\hdots,x_{n+1}]/(q-1),\quad \Delta:=(-1)^{\frac{n(n+1)}{2}}\prod_{i=1}^{n+1} a_i.
    \end{gathered}
\]
By the adjunction formula we have $\omega_{Q/k}\cong \struct(-n)$ and $\omega_{Q^\infty/k}\cong \struct(-n+1)$, where we denote by $\struct(-1)$ the restriction of the tautological bundle from $\PP^{n+1}$ to the corresponding quadric. Furthermore, $\omega_{Q^o/k}\cong \struct_{Q^o}$ using the splitting~\eqref{eq:tangent_stably_trivial}, and we choose the isomorphism to be given by the globally defined nowhere vanishing top dimension differential form 
\begin{equation}\label{eq:theta_form}
\theta:=\sum_{i=1}^{n+1}(-1)^{i-1} x_idx_1\wedge dx_2\wedge \hdots\wedge\widehat{dx_i}\wedge\hdots \wedge dx_{n+1}.
\end{equation}

\begin{theorem} \label{thm:n_odd_necessary}
    Suppose that $n\ge 3$ is odd, $\Delta \notin (k^\times)^2$, and  $T_{Q^o}$ is trivial.
    Then $-a_1\in [D(q)^2]$.
\end{theorem}
\begin{proof} 
    Consider the following commutative diagram.
    \[
    \xymatrix{
    & \K_1(Q)\ar[r]^{j_K^*}\ar[d]^{c_{1,n+1}} & \K_1(Q^o)\ar[d]^{c_{1,n+1}}  & \\
    \H^{n-1}(Q^\infty,\K^\M_{n}) \ar[r]^{i_*} \ar[rd]_{\pi^\infty_*} & \H^n(Q,\K^\M_{n+1}) \ar[r]  \ar[d]^{\pi_*} & \H^n(Q^0,\K^\M_{n+1})\ar[r] & 0 \\
    & k^\times & &
    }
    \]
    Here $j_K^*$ is the pullback along the open embedding, the middle row is the localization exact sequence, $c_{1,n+1}$ are higher Chern classes, and $\pi_*$, $\pi^\infty_*$ are pushforwards to $\Spec k$. Consider the cycle $\alpha\in\H^n(Q^0,\K^\M_{n+1})$ corresponding to the unimodular row $(x_1,x_2,\hdots,x_{n+1})$,
    \[
    \alpha:=\bar{\phi}(x_1,x_2,\hdots,x_{n+1}) = \{t\}_Z,\quad
    Z:=\{x_2=\hdots=x_{n+1}=0\}\cong \Spec k[t]/(a_1t^2-1)\subseteq Q^o.
    \]
    The exact sequence~\eqref{eq:stably_trivial_as_cokernel} and triviality of $T_{Q^o}$ yield that $\alpha \in \image c_{1,n+1}$. If $\alpha'$ is the tautological lift of $\alpha$ to $Q$, then Lemma~\ref{lem:K-theory_quadrics_surjection} yields $\alpha' \in \image i_* + c_{1,n+1}(\K_1(Q))$.
    Hence 
    \[
    -a_1^{-1}=\pi_*(\alpha') \in \image \pi_*^\infty + \pi_* c_{1,n+1}(\K_1(Q)) = [D(q)^2],
    \]
    where for the last equality we used Lemmas~\ref{lem:image_transfer_KMW1} and~\ref{lem:K-theory_quadrics_push_higher}. The claim follows.
\end{proof}

\begin{theorem} \label{thm:n_odd_sufficient}
    Suppose that $n\ge 3$ is odd, the pushforward $\pi_*\colon \H^n(Q,\K^\M_{n+1})\to k^\times$ is injective and $-a_1\in [D(q)^2]$.
    Then $T_{Q^o}$ is trivial.
\end{theorem}
\begin{proof}
    In view of Lemma~\ref{lem:completable_trivial} and bijection~\eqref{eq:obstruction_odd_triviality} it suffices to show that for the unimodular row $(x_1,x_2,\hdots,x_{n+1})$ one has $\bar{\phi}(x_1,x_2,\hdots,x_{n+1})=0$ in $\H^n(Q^0,\K^\M_{n+1})$.

    As above, we have
    \[
    \bar{\phi}(x_1,x_2,\hdots,x_{n+1}) = \{t\}_Z,\quad
    Z:=\{x_2=\hdots=x_{n+1}=0\}\cong \Spec k[t]/(a_1t^2-1)\subseteq Q^o.
    \]
    Consider the following commutative diagram
    \[
    \xymatrix{
     \H^{n-1}(Q^\infty,\K^\M_{n}) \ar[r]^{i_*} \ar[rd]_{\pi^\infty_*} & \H^n(Q,\K^\M_{n+1}) \ar[r]^{j^*}  \ar[d]^{\pi_*} & \H^n(Q^0,\K^\M_{n+1})\ar[r] & 0. \\
    & k^\times & &
    }
    \]
    Here the horizontal row is the localization exact sequence and $\pi_*$, $\pi^\infty_*$ are pushforwards to $\Spec k$. The cycle $\{t\}_Z$ tautologically lifts to $\H^n(Q,\K^\M_{n+1})$. We have $\pi_*(\{t\}_Z)=-a_1^{-1}$, and $-a_1^{-1}\in \image \pi^\infty_*$ by Lemma~\ref{lem:image_transfer_KMW1}. Then injectivity of the pushforward yields that $\{t\}_Z \in \image i_*$, hence $j^*(\{t\}_Z)=0$ and the claim follows.
\end{proof}

\begin{remark} \label{rem:injectivity}
    The injectivity condition of Theorem~\ref{thm:n_odd_sufficient} is satisfied for $\tilde{q}=q_1\perp q_2$ with $q_1$ being a Pfister neighbor and $\dim q_2\le 2$ by \cite[Corollary~6.3]{CM01} (see also \cite[Proposition~2]{Ros88}). In particular, it holds for all smooth projective quadrics of dimension $3$.
\end{remark}

\begin{corollary}\label{cor:dim3}
    Suppose $n=3$. Then the tangent bundle $T_{Q^o}$ is trivial if and only if one of the following holds: (i) $\Delta\in (k^\times)^2$,
        (ii) $-a_1\in [D(q)^2]$.
\end{corollary}
\begin{proof}
    One direction follows from Example~\ref{ex:dim3}, Theorem~\ref{thm:n_odd_sufficient} and Remark~\ref{rem:injectivity}. The other direction follows from Theorem~\ref{thm:n_odd_necessary}.
\end{proof}

\begin{recollection}
    Let $X\in \Sm_k$ and $\mc{V}$ be a rank $n$ vector bundle over $X$. There is a natural Euler class
    \[
    e(\mc{V})\in \H^n(X,\K^\MW_n\otimes(\det \mc{V})^\vee),
    \]
    see e.g. \cite[Chapter~8]{Mor12} and \cite[\S3]{AF16}.
    Suppose that $s$ is a section of $\mc{V}$ transverse to the zero section and $i\colon Z\to X$ is the inclusion of the zero locus, then by \cite[Corollary~5.21]{BW23} we have $e(\mc{V})=i_*(1)$ for the pushforward
    \[
    i_*\colon  \H^{0}(Z,\K^\MW_0)\cong \H^{0}(Z,\K^\MW_0\otimes \omega_{Z/k}\otimes(\omega_{X/k})_Z^{\vee}\otimes (\det \mc{V})_Z^\vee)\to  \H^{n}(X,\K^\MW_n\otimes (\det \mc{V})^\vee),
    \]
   where the isomorphism is induced by $\mc{N}_{Z/X}\cong \mc{V}_Z$ given by $s$. In particular, if $n=\dim X$ then 
    \[
    e(\mc{V}) = \langle 1\rangle_Z \otimes \lambda_Z,
    \]
    with $\lambda_Z\in \Gamma(Z,\omega^\vee_{X/k} \otimes \det \mc{V}^\vee)$
    being the nowhere vanishing section arising from the determinant of the isomorphism $(ds)_Z\colon (T_X)|_Z\xrightarrow{\simeq}\mc{V}|_Z$, see \cite[Example~5.6]{Lev20} for more details.
    
    If $X$ is affine, $n=\dim X$ and $\det \mc{V}\cong \struct_X$ then $\mc{V}$ has a nowhere vanishing section if and only if $e(\mc{V})=0$ by \cite[Chapter~8]{Mor12}, \cite[Theorem~1]{AHW17} and \cite[Theorem~5.6]{AF16}.
\end{recollection}

The following result can be obtained combining \cite[\S 3]{DS15} and \cite[proof of Theorem 3.2.6]{AF22}, but we reprove it for the sake of completeness.

\begin{lemma}\label{lem:etaEulerclass}
    Suppose $n\ge 2$ is even. Then $e(T_{Q^o})=\eta \phi (x_1,x_2,\hdots,x_{n+1}) $ under the isomorphism $\H^n(Q^o,\K^\MW_n) \cong  \H^n(Q^o,\K^\MW_n \otimes \omega_{Q^o/k})$ induced by $\theta$ from~\eqref{eq:theta_form}.
\end{lemma}
\begin{proof}
    We have 
    \begin{gather*}
    \phi(x_1,x_2,\hdots,x_{n+1}) = [t]_Z\otimes \gamma_Z^\vee\in  \H^n(Q^o,\K^\MW_{n+1}),\quad \gamma:={dx}_2\wedge {dx}_3\wedge\hdots\wedge {dx}_{n+1},\\
    Z:=\{x_2=\hdots=x_{n+1}=0\}\cong \Spec k[t]/(a_1t^2-1)\subseteq Q^o.
    \end{gather*}
    Multiplying by $\eta$, we obtain $\eta\phi(x_1,x_2,\hdots,x_{n+1})= (\langle t\rangle-1)_Z\otimes \gamma_Z^\vee$. Now, $1_Z\otimes \gamma_Z^\vee$ is the cycle associated with the complete intersection $(x_2,x_3,\hdots,x_{n+1})$ and thus vanishes by \cite[Th\'eor\`eme 10.3.5]{Fas08}. Consequently,
    \begin{equation}\label{eqn:etaphi}
    \eta\phi(x_1,x_2,\hdots,x_{n+1})=\langle t\rangle_Z\otimes \gamma_Z^\vee.
    \end{equation} 
Pick the section $s$ of $T_{Q^o}$ given by the row $(1-a_1x_1^2,-a_1x_1x_2,\hdots,-a_1x_1x_{n+1})$. Its zero locus is precisely $Z$, and the intersection of $s(Q^o)$ with the zero section is transversal. We have
    \[
    e(T_{Q^o})= \langle 1\rangle_Z \otimes \left((-a_1t)^{-n} \gamma^\vee\otimes \gamma\right) = \langle (-a_1t)^{-n}\rangle_Z=\langle 1\rangle_Z,
    \]
    the last equality uses that $n$ is even. It is now straightforward to check that the image of $\langle t\rangle_Z\otimes \gamma_Z^\vee$ under the isomorphism given by $\theta$ is $\langle 1\rangle_Z$. 
\end{proof}

\begin{lemma} \label{lem:conic_inj}
    Let $Q=\{q=0\}\subseteq \PP^2$ be a smooth conic over $k$. Then the pushforward homomorphism
    \[
    \pi_*\colon \H^1(Q,\K^\MW_2\otimes \omega_{Q/k}) \to \K_1^\MW(k)
    \]
    is injective.
\end{lemma}

\begin{proof}
We know from Corollary \ref{cor:K1symplecticinjective}  that the pushforward homomorphism
\[
p_*\colon \mathrm{GW}_1^{[2]}(Q,\omega_{Q/k})\to \mathrm{GW}^{[1]}_1(k)
\]  
is injective. On the other hand, we know from \cite[\S 4]{AF17} that there is a canonical morphism of sheaves $\KMW_i\to \mathrm{GW}^{[i]}_i$ that restricts to an isomorphism $\KMW_i(F)\to \mathrm{GW}^{[i]}_i(F)$ for any field $F$ and any $i\leq 3$. To identify $\mathrm{GW}_1^{[2]}(Q,\omega_{Q/k})$, we use the coniveau spectral sequence for Grothendieck--Witt groups (\cite[\S 3]{FS09}). Since $\mathrm{GW}_1^{[2]}(k(Q))=0$, it follows that $\mathrm{GW}_1^{[2]}(Q,\omega_{Q/k})$ is isomorphic to the cokernel of the residue map
\[
\mathrm{GW}_2^{[2]}(k(Q),\omega_{k(Q)/k})\to \bigoplus_{x\in Q_{(0)}}\mathrm{GW}_1^{[1]}(k(x))
\]
which coincides with the group $ \H^1(Q,\K^\MW_2\otimes \omega_{Q/k})$ by \cite[\S 3.4]{FS09}.    
\end{proof}

\begin{definition}
Recall that for every $m$ we have the following short exact sequences of sheaves.
\begin{gather*}
 0 \to \I^{m+1} \to \KMW_m\to \KM_m\to 0\\
 0 \to 2\KM_{m} \to \KMW_m\to \I^m\to 0\\
    0 \to {}_2\KM_{m} \to \KM_{m}\to 2\KM_{m}\to 0
\end{gather*}
which give rise to the long exact sequences of (twisted) cohomology groups. We denote by $\mc{F}\colon\KMW_m\to \KM_m$ the projection and by $\mc{H}\colon \KM_{m}\to 2\KM_{m}\to \KMW_m$ the composition, to which we refer as \textit{hyperbolization}. For a field $L$ and $\{u\}\in \KM_1(L)$ corresponding to $u\in L^\times$ one has $\mc{H}(\{u\})=[u^2]=h[u]$, where $h=\langle 1,-1\rangle$.
\end{definition}

\begin{proposition}\label{prop:hyperbolization_zero}
    Let $n\ge 1$. Then $\mc{H}\colon \H^n(Q^o,\K^\M_{n+1}) \to \H^n(Q^o,\K^\MW_{n+1})$
    is zero.
\end{proposition}
\begin{proof}
    Consider the following commutative diagram.
    \[
    \xymatrix{
    &\H^n(Q,\K^\M_{n+1}) \ar[r]^{j^*} \ar[d]_{\mc{H}} & \H^n(Q^o,\K^\M_{n+1}) \ar[d]^{\mc{H}} \\
    \H^{n-1}(Q^\infty,\K^\MW_{n}\otimes \omega_{Q^\infty/k}) \ar[r]^(0.55){i_*} & \H^n(Q,\K^\MW_{n+1}\otimes \omega_{Q/k}) \ar[r] & \H^n(Q^o,\K^\MW_{n+1})
    }
    \]
    Since $j^*$ is surjective it suffices to show that for $z\in Q^o$, $u\in k(z)^\times$ one has $\mc{H}(\{u\}_z)\in \image i_*$ for the respective cycle $\{u\}_z\in \H^n(Q,\K^\M_{n+1})$. 

    First suppose $\tilde{q}$ to  be isotropic. Then pushforward homomorphisms $\tilde{\pi}_*\colon \H^n(Q,\K^\MW_{n+1}\otimes \omega_{Q/k}) \to \K^\MW_1(k)$ and ${\pi}_*\colon \H^n(Q,\K^\M_{n+1}) \to \K^\M_1(k)$ are isomorphisms by Corollary~\ref{cor:stable_birational_invariant}, so it suffices to show that for every $u\in k^\times$ one has $h[u]\in \image \pi^\infty_*$ for the morphism $\pi^\infty\colon Q^\infty\to \Spec k$. If $q$ is isotropic then $\pi^\infty_*$ is surjective (and an isomorphism if $n\ge 2$) and the claim is trivial. Otherwise pick a degree 2 point $z\in Q^\infty$, then $k(z)\cong k[\sqrt{t}]$ for some $t\in k^\times$. We have $\pi^\infty_*(\langle \sqrt{t}\rangle_z[u]_z)=\pi^\infty_*(\langle \sqrt{t}\rangle_z)[u] =h[u]$, proving the claim.
    
    Now suppose $\tilde{q}$ to  be anisotropic. In view of Lemma~\ref{lem:deg2_generation} we may assume $[k(z):k]=2$. Then there is a dimension $2$ subspace $W\subseteq k^{n+2}$ such that $z\in \mathbb{P}(W)\subseteq \PP^{n+1}$. Choose an arbitrary $W\subseteq V\subseteq k^{n+2}$ of dimension $3$. Since $\tilde{q}$ is anisotropic then $\tilde{q}|_V$ is non-degenerate and $C:=\PP(V)\cap Q$ is a conic. We have $Q^\infty \cap C = \{x_0=0\} \cap C  =\{w\}$ with $w\cong \Spec k[\sqrt{t}]$ for some $t\in k^\times$. Consider the cycle $\langle \sqrt{t} \rangle_w [N_{k(z)/k}(u)]_w \in \H^{1}(C,\K^\MW_{2}\otimes \omega_{C/k})$. Then 
    \begin{multline*}
    (\pi_C)_*(\langle \sqrt{t} \rangle_w [N_{k(z)/k}(u)]_w) = (\pi_C)_*(\langle \sqrt{t} \rangle_w) [N_{k(z)/k}(u)] =h [N_{k(z)/k}(u)] = \\
    =\mc{H}((\pi_C)_*(\{u\}_z)) = (\pi_C)_*(\mc{H}(\{u\}_z))
    \end{multline*}
    for the projection $\pi_C\colon C\to \Spec k$.
    Lemma~\ref{lem:conic_inj} yields $\langle \sqrt{t} \rangle [N_{k(z)/k}(u)]_w = \mc{H}(\{u\}_z)$ in $\H^{1}(C,\K^\MW_{2}\otimes \omega_{C/k})$. Pushing this forward along the embedding $C\to Q$ we obtain that $\mc{H}(\{u\}_z)$ equals in $\H^n(Q,\K^\MW_{n+1}\otimes \omega_{Q/k})$ to the cycle $\langle \sqrt{t} \rangle [N_{k(z)/k}(u)]_w$ supported on $Q^\infty$, proving the claim.
\end{proof} 

\begin{corollary} \label{cor:inj_proj_I}
    Let $n\ge 1$. Then the homomorphism $\H^n(Q^o,\K^\MW_{n+1}) \to \H^n(Q^o,\I^{n+1})$ is an isomorphism.
\end{corollary}
\begin{proof}
    Consider the exact sequence
    \[
    \H^n(Q^o,\K^\M_{n+1}) \xrightarrow{f} \H^n(Q^o,2\K^\M_{n+1}) \to \H^{n+1}(Q^o,{}_2\K^\M_{n+1}).
    \]
    The last group vanishes by dimension reasons, so $f$ is surjective. Hyperbolization $\mc{H}$ factors as
    \[
    \H^n(Q^o,\K^\M_{n+1}) \xrightarrow{f} \H^n(Q^o,2\K^\M_{n+1})\xrightarrow{\mc{H}'} \H^n(Q^o,\K^\MW_{n+1})
    \]
    hence Proposition~\ref{prop:hyperbolization_zero} yields that $\mc{H}'$ is zero. The claim follows from the exact sequence
    \[
    \H^n(Q^o,2\K^\M_{n+1})\xrightarrow{\mc{H}'} \H^n(Q^o,\K^\MW_{n+1}) \to \H^n(Q^o,\I^{n+1})\to 0. \qedhere
    \]
\end{proof}

\begin{proposition} \label{prop:inj_inj_I}
    The homomorphism $\H^n(Q^o,\I^{n+1}) \to \H^n(Q^o,\K^\MW_{n})$ is injective.
\end{proposition}
\begin{proof}
    If $n=0$ the claim is clear, while for $n=1$ the claim follows from the long exact sequence associated to the short exact sequence $\I^{2} \to \KMW_1\to \KM_1$ of sheaves and an obvious surjectivity of
    \[
    \H^0(Q^o,\K^\MW_{1})\to \H^0(Q^o,\K^\M_{1})\cong k[Q^o]^\times,
    \]
    so we may assume $n\ge 2$.

    Consider the localization sequence
    \[
    \H^{n-1}(Q,\K^\M_{n}) \xrightarrow{j_K^*} \H^{n-1}(Q^o,\K^\M_{n}) \to  \H^{n-1}(Q^\infty,\K^\M_{n-1}) \xrightarrow{i_*}  \H^{n}(Q,\K^\M_{n})
    \]
    The homomorphism $i_*$ is injective by \cite[Corollary~70.4, Theorem~71.2]{EKM08}, hence $j_K^*$ is surjective. Consider the following commutative diagram.
    \[
    \xymatrix{
     \H^{n-1}(Q,\K^\M_{n}) \ar[d]_{j^*_K} \ar[r] & \H^{n-1}(Q,\kM_{n}) \ar[r]^(0.43){\delta} & \H^{n}(Q,\I^{n+1}\otimes \omega_{Q/k}) \ar[d]^{j^*_\I} \\
    \H^{n-1}(Q^o,\K^\M_{n}) \ar[rr]^{\partial_o} & & \H^{n}(Q^o,\I^{n+1}) \ar[r] & \H^{n}(Q^o,\KMW_{n})
    }
    \]
    Here the lower row is exact, and the unnamed morphism in the first row is the mod $2$ reduction. The claim of the proposition is equivalent to $\partial_o=0$. Since $j^*_K$ is surjective by the above, it suffices to show that $j_\I^*\circ\delta=0$. By Corollary~\ref{cor:rost} this is equivalent to the following two vanishings.
        \begin{enumerate}
        \item $j^*_{I}   \delta (\{u\}\cdot [h]^{n-1})=0$ for $u\in k^\times$ and the class of the hyperplane section $[h]\in\Ch^{1}(Q)$,
        \item $(\pi_{L/k})_*j^*_{I}   \delta  (\{u\}\cdot [\ell_L])=0$ for $u\in L^\times$, the class of a line $[\ell_L]\in\Ch^{n-1}(Q_L)$ and $L/k$ finite field extension such that $i(\tilde{q}_L)\ge 2$.
    \end{enumerate}
    Carrying $\{u\}$ as $\langle\langle u\rangle\rangle$ out in the above formulas it suffices to show $j^*_{I}   \delta ([h]^{n-1})=0$ and $j^*_{I}   \delta  ([\ell_L])=0$. For this consider the following commutative diagram.
     \[
    \xymatrix{
    & \H^{n-1}(Q,\I^{n-1} \otimes \omega_{Q/k}) \ar[d]_{\mc{F}} & \\
    \Ch^{n-2}(Q^\infty) \ar[r]^{i_*} &\Ch^{n-1}(Q) \ar[r]^{j_K^*} \ar[d]_{\delta} & \Ch^{n-1}(Q^o) \ar[d]^{\delta_o} \\
    &\H^{n}(Q,\I^{n} \otimes \omega_{Q/k}) \ar[r]^{j_{I}^*} & \H^{n}(Q^o,\I^{n})
    }
    \]
    Here the column and the row are exact, and we denote by the same letters as before the morphisms one weight below.
    
    For the first vanishing we have
    \[
    j^*_{I} \delta ([h]^{n-1})= \delta_o j^*_{K} ([h]^{n-1}) = \delta_o j^*_{K} i_*([h]^{n-2})= 0.
    \]
    
    For the second vanishing we base change to $L$ and pick a linear embedding $s\colon \PP^1_L\to Q_L$. Consider the following commutative diagram.
    \[
    \xymatrix{
    \H^{0}(\PP^1_L,\mr{W}) \ar[r]_(0.4){\cong}^(0.4){\theta} \ar[dr]_{\mc{F}_{\PP^1}}  &\H^{0}(\PP^1_L,\mr{W} \otimes \omega_{\PP^1_L/k}) \ar[d]^{\mc{F}_{\PP^1,\omega}} \ar[r]^(0.46){\tilde{s}_*} &\H^{n-1}(Q_L,\I^{n-1} \otimes \omega_{Q_L/k}) \ar[d]^{\mc{F}}\\
    &\Ch^{0}(\PP^1_L) \ar[r]^{s_*} &\Ch^{n-1}(Q_L)
    }
    \]
    Here the vertical morphisms are induced by the morphism of sheaves $\I^*\to \overline{\I}{}^*\cong\kM_*$ and isomorphism $\theta$ is induced by the choice of the square root $\omega_{\PP^1_L/k}\cong \struct(-2)\cong \struct(-1)^{\otimes 2}$. We have $\mc{F}_{\PP^1}(1)=1$ hence $s_*(1) =  \mc{F}\tilde{s}_* \theta (1)$ and thus
    \[
    j^*_{I} \delta [\ell_L]= j^*_{I} \delta s_*(1) =j^*_{I} (\delta \mc{F}) \tilde{s}_* \theta (1) =0. \qedhere
    \]
\end{proof}

\begin{theorem} \label{thm:even_parallel}
    Suppose $n>0$ is even. Then the following are equivalent.
    \begin{enumerate}
        \item $T_{Q^o}$ has a nowhere vanishing section,
        \item$T_{Q^o}$ is trivial.
    \end{enumerate}
    If $\pi_*\colon \CHW_0(Q)\to \GW(k)$ is injective, then the above are also equivalent to
    \begin{enumerate}[resume]
        \item $-a_1 \in [D(q)^2]$.
    \end{enumerate}
\end{theorem}
\begin{proof}
    The implication $(2)\Rightarrow (1)$ is trivial. For $(1)\Rightarrow (2)$ suppose that $T_{Q^o}$ has a nowhere vanishing section. Hence $e(T_{Q^o})=0$. Multiplication by $\eta$ factors as
    \[
    \H^n(Q^o,\K^\MW_{n+1}) \to\H^n(Q^o,\I^{n+1}) \to \H^n(Q^o,\K^\MW_{n})
    \]
    which is injective by Corollary~\ref{cor:inj_proj_I} and Proposition~\ref{prop:inj_inj_I}. Since $e(T_{Q^o})=\eta \phi(x_1,x_2,\hdots,x_{n+1})$ by Lemma~\ref{lem:etaEulerclass}, we have $\phi(x_1,x_2,\hdots,x_{n+1})=0$ and $(2)$ follows from the bijection~\eqref{eq:obstruction_general_triviality}.

    The equivalence $(1)\Leftrightarrow (3)$ under the assumption of injectivity of $\pi_*$ follows by the same argument as in the proof of \cite[Theorem~4.1]{AL25}. Note that since $n$ is even, $-a_1 \in [D(q)^2]$ if and only if $-\prod_{i=1}^{n+1}a_i \in [D(q)^2]$.
\end{proof}

\begin{remark}
    We know that the pushforward homomorphism $\pi_*$ from Theorem~\ref{thm:even_parallel} is injective if $Q^o(k)\neq \emptyset$ by Corollary~\ref{cor:stable_birational_invariant} and also if $n\le 3$ by Corollary~\ref{cor:explicitkernel}. Furthermore, if $\pi_*$ is injective, Theorem~\ref{thm:even_parallel} admits a much more straightforward proof similar to the one given for Theorem~\ref{thm:n_odd_sufficient}. 
\end{remark}

\begin{remark}
     In case $Q^o(k)\neq \emptyset$ we have $[D(q)^2]=[D(q)]$ and $(3)$ of Theorem~\ref{thm:even_parallel} is equivalent to $-1 \in [D(q)]$ as in \cite[Theorem~4.1.(3)]{AL25}.
\end{remark}

\section{Triviality in higher dimension: vanishing birational invariants}\label{sec: vanishingbirational}

In this section we construct a family of non-rational affine quadrics with trivial tangent bundles. This is achieved by establishing the vanishing of groups of zero cycles with coefficients in all homotopy modules using some explicit birational geometry of quadrics.
 
\begin{definition}
    The \textit{Witt index} $i(q)$ of a regular quadratic form on a vector space $V$ over $k$ is the dimension of a maximal totally isotropic subspace of $V$.
    The \textit{first higher Witt index} $i_1(q)$ of an anisotropic quadratic form $q$ of $\dim q\ge 2$ over $k$ is the Witt index $i(q_{k(Q)})$ of the quadratic form $q_{k(Q)}$ over the function field $k(Q)$ of the projective quadric $\{q=0\}$. Equivalently, $i_1(q)$ is the minimal positive Witt index $i(q_L)$ for all field extensions $L/k$.
\end{definition}

\begin{example}
    Let $q$ be a regular anisotropic quadratic form over $k$. Suppose that $q$ is a Pfister neighbor of dimension $2^{n-1}+m$, $1\le m\le 2^{n-1}$. Then $i_1(q)= m$ by \cite[\S 3, Corollary 2]{Hof95}.
\end{example}

\begin{theorem} \label{thm:birational_vanishing}
    Let $q$ be a regular anisotropic quadratic form of dimension $n+2$ over $k$ and consider $Q:=\{q=0\}\subseteq \PP^{n+1}$, $Q^\infty:=Q\cap \{x_{0}=0\}$, $Q^o:=Q\setminus Q^\infty$.
    Suppose $i_1(q)\ge 2$. Then $A_0(Q^o,M_*)=0$ for every homotopy module $M_*$. 
\end{theorem}
\begin{proof}
Let $i\colon Q^\infty\to Q$ be the embedding and consider the localization exact sequence
\[
A_0(Q^\infty,M_*) \xrightarrow{i_*} A_0(Q ,M_*) \to A_0(Q^o,M_*) \to 0.
\]
The assumption on the Witt index yields that $Q^\infty$ is stably birational to $Q$, see e.g. \cite[\S 1]{Tot09}, so the claim basically follows since $A_0(-,M_*)$ is a stable birational invariant by Corollary~\ref{cor:stable_birational_invariant}.
Below are the details showing that $i_*$ is indeed an isomorphism.

The vector space $k(Q)^{n+2}$ admits a $2$-dimensional subspace isotropic with respect to the quadratic form $q$ by the assumption on the Witt index. Hence the hyperplane $\{x_{0}=0\}$ admits a $1$-dimensional isotropic space, i.e. $Q^\infty_{k(Q)}$ admits a rational point. It follows that there exists a non-empty open subset $U\subseteq Q$ and a section $  U\to Q^\infty \times U$ for the projection $Q^\infty \times U\to U$. The usual stereographic projection applied in the family over $U$ gives rise to a birational equivalence  $Q^\infty \times U \dashrightarrow \PP^{n-1}\times U$ over $U$. Using the embedding $U\subseteq Q$ we obtain a commutative diagram
\[
\xymatrix{
Q^\infty \times Q \ar[dr]_{\pi_Q} \ar@{-->}[rr]^{\Phi} & & \PP^{n-1}\times Q \ar[dl]^{\pi'_Q} \\
& Q &
}
\]
with $\Phi$ being a birational equivalence and $\pi_Q, \pi'_Q$ being the projections. In the remaining reasoning we use the notation of \cite[2.2.6 and~2.2.8]{Fel22} and implicitly invoke \cite[Corollary~2.2.11]{Fel22}. Since the diagram commutes we have
\[
[\pi_Q]_*=[\pi'_Q]_* \circ [\Phi]_* \colon A_0(Q^\infty \times Q,M_*) \to A_0(Q,M_*).
\]
Furthermore, $[\Phi]_*$ is an isomorphism since $\Phi$ is a birational equivalence, while $[\pi'_Q]_*=(\pi'_Q)_*$ is an isomorphism by Lemma~\ref{lem:push_proj_iso}, yielding that $[\pi_Q]_*$ is an isomorphism as well. Since $Q_{k(Q^\infty)}$ tautologically has a rational point then a similar reasoning yields that for the projection $\pi_{Q^\infty}\colon Q^\infty \times Q \to Q^\infty$ the homomorphism
\[
[\pi_{Q^\infty}]_*\colon A_0(Q^\infty \times Q,M_*) \to A_0(Q^\infty,M_*)
\]
is an isomorphism as well. The morphism $s:=(\id,i)\colon Q^\infty \to Q^\infty\times Q$ is a section to $\pi_{Q^\infty}$ hence $[s]_*$ is also an isomorphism. Finally, $i=\pi_Q\circ s$ hence $i_*=[i]_*$ is an isomorphism as well.
\end{proof}

\begin{theorem}\label{thm:vanishing_birational}
    Let $q=\sum_{i=1}^{n+1}a_ix_i^2$ be a regular quadratic form over $k$, put $\tilde{q}:=-x_0^2+q$ and $ Q^o:=\{q=1\}$. Suppose that $\tilde{q}$ is anisotropic and $i_1(\tilde{q})\ge 2$.
    Then every rank $n$ stably trivial vector bundle over $Q^o$ is trivial. In particular, the tangent bundle $T_{Q^o}$ is trivial.
\end{theorem}
\begin{proof}
    Since $\omega_{Q^0/k}$ is trivial, Theorem~\ref{thm:birational_vanishing} yields vanishing
    \[
    \mathrm{H}^n(Q^o,\mathrm{K}_{n+1}^{\mr{MW}}) \cong A_0(Q^o,\mathrm{K}_{1}^{\mr{MW}}) =0.
    \]
    Hence bijection~\eqref{eq:obstruction_general_triviality} yields that $\Um_{n+1}(R)/\GL_{n+1}(R)\cong \Vect^o_n(Q^o)$
    is trivial. Serre's cancellation then yields that all rank $n$ stably trivial vector bundles over $Q^o$ are trivial.
\end{proof}

\section{A useful theorem}

The object of this section is the following theorem. 

\begin{theorem}\label{thm:pushforwardequivalence}
Let $p\colon Q\to \Spec(k)$ be a smooth projective quadric of dimension $n\geq 1$. Let $p_*^{\mathrm{GW}}\colon \CHW_0(Q)\to \CHW_0(k)$ and $p_*^{\mathrm{W}}\colon \H^n(Q,\mathrm{W}\otimes \omega_{Q/k})\to \mathrm{W}(k)$ be the pushforward maps. Then, the homomorphism $\CHW_0(Q)\to \H^n(Q,\mathrm{W}\otimes \omega_{Q/k})$ induces an isomorphism $\ker(p_*^{\mathrm{GW}})\simeq \ker(p_*^{\mathrm{W}})$.
\end{theorem}

The proof of the theorem will consist in a few computations in cohomology that we now start. We first note that if $Q$ is isotropic both maps are isomorphisms by Corollary~\ref{cor:stable_birational_invariant} and there is nothing to prove. We thus assume from now on that $Q$ is anisotropic.

The exact sequence of sheaves 
\[
0\to \mathrm{I}^n\otimes\omega_{Q/k}\to \mathrm{I}^{n-1}\otimes\omega_{Q/k}\to \overline{\mathrm{I}}{}^{n-1}\to 0
\]
induces an exact sequence of the form
\[
\H^{n-1}(Q,\mathrm{I}^{n-1}\otimes\omega_{Q/k})\to \H^{n-1}(Q,\overline{\mathrm{I}}{}^{n-1})\xrightarrow{\partial} \H^n(Q,\mathrm{I}^n\otimes \omega_{Q/k})\to \H^n(Q,\mathrm{I}^{n-1}\otimes \omega_{Q/k})
\]
Recall from \cite[Theorem 1.1]{Tot09} that $\H^{n-1}(Q,\overline{\mathrm{I}}{}^{n-1})=\mathrm{Ch}^{n-1}(Q)$.
\begin{lemma}
Let $[h]\in \Ch^1(Q)$ be the class of a hyperplane section. Then, the image of $\partial$ is generated by $\partial([h]^{n-1})$.
\end{lemma}

\begin{proof}
We know from Lemma \ref{lem:baby_rost} that $\mathrm{Ch}^{n-1}(Q)$ is generated by $[h]^{n-1}$, together with the images under pushforward of lines $[\ell_1]\in \Ch^{n-1}(Q_L)$ where $L/k$ is such that $q_L=q'\perp 2\mathbb{H}$ for some quadratic form $q'$ over $L$. In particular $q_L$ is isotropic and the pushforward homomorphisms 
\[
p_*\colon \H^n(Q_{L},\mathrm{I}^n\otimes \omega_{Q_L/k})\to \mathrm{W}(L),\phantom{i}p_*\colon \H^n(Q_{L},\mathrm{I}^{n-1}\otimes \omega_{Q_L/k})\to \mathrm{W}(L)
\]
are both isomorphisms. Thus, the map $\H^n(Q_{L},\mathrm{I}^n\otimes \omega_{Q_L/k})\to \H^n(Q_{L},\mathrm{I}^{n-1}\otimes \omega_{Q_L/k})$ is an isomorphism and the connecting homomorphism $\partial_L\colon  \H^{n-1}(Q_{L},\overline{\mathrm{I}}{}^{n-1})\to \H^n(Q_{L},\mathrm{I}^n\otimes \omega_{Q_L/k})$ is trivial. The commutative diagram
\[
\xymatrix{ \H^{n-1}(Q_{L},\overline{\mathrm{I}}^{n-1})\ar[r]^-{\partial_L}\ar[d] & \H^n(Q_{L},\mathrm{I}^n\otimes \omega_{Q_L/k})\ar[d] \\
\H^{n-1}(Q,\overline{\mathrm{I}}^{n-1})\ar[r]_-{\partial} & \H^n(Q,\mathrm{I}^n\otimes \omega_{Q/k}), }
\]
where the vertical maps are transfer homomorphisms, allows to conclude that the images of the lines have trivial image under $\partial$. The claim follows.
\end{proof}

\begin{corollary}
The obvious morphism $\mathrm{I}^n\otimes\omega_{Q/k}\to \overline{\mathrm{I}}{}^n \oplus \mathrm{I}^{n-1}\otimes\omega_{Q/k}$ induces an isomorphism
\[
\H^n(Q,\mathrm{I}^n\otimes\omega_{Q/k})\simeq  \H^n(Q, \overline{\mathrm{I}}{}^n) \oplus \H^n(Q,\mathrm{I}^{n-1}\otimes\omega_{Q/k}).
\]
\end{corollary}

\begin{proof}
By the above lemma, we have an exact sequence
\[
 \Z/2\cdot \partial([h]^{n-1})\to \H^n(Q,\mathrm{I}^n\otimes\omega_{Q/k})\to \H^n(Q,\mathrm{I}^{n-1}\otimes\omega_{Q/k})\to \H^n(Q,\overline{\mathrm{I}}{}^{n-1})=0
\]
and it suffices to prove that the composite 
\[
 \Z/2\cdot \partial([h]^{n-1})\to \H^n(Q,\mathrm{I}^n\otimes\omega_{Q/k})\to \H^n(Q, \overline{\mathrm{I}}{}^n)=\mathrm{Ch}^n(Q)=\Z/2\Z
\]
is an isomorphism to conclude. To see this, we first deduce from \cite[Theorem 3.4.1]{AF15} and \cite[Theorem 1.1]{Tot03} that the image of $\partial([h]^{n-1})$ under the composite is equal to $\mathrm{Sq}^2_{\omega}([h]^{n-1})$, where the latter is the twisted Steenrod operator with $\omega:=\omega_{Q/k}$. 

To compute $\mathrm{Sq}^2_{\omega}([h]^{n-1})$, we first consider the operation $\mathrm{Sq}^2_{\omega'}\colon \mathrm{Ch}^{n-1}(\mathbb{P}^{n+1})\to   \mathrm{Ch}^{n}(\mathbb{P}^{n+1})$ where $\omega'=\omega_{\mathbb{P}^{n+1}/k}=\mathcal{O}(-n-2)$, which a direct computation shows to be an isomorphism. If $i\colon Q\to \mathbb{P}^{n+1}$ is the inclusion, we may use the commutative diagram 
\[
\xymatrix{\mathrm{Ch}^{n-1}(\mathbb{P}^{n+1})\ar[r]^-{\mathrm{Sq}^2_{\omega'}}\ar[d]_-{i^*} &\mathrm{Ch}^{n}(\mathbb{P}^{n+1})\ar[d]^-{i^*}  \\
\mathrm{Ch}^{n-1}(Q)\ar[r]_-{\mathrm{Sq}^2_{\omega}} &\mathrm{Ch}^{n}(Q)}
\] 
and the fact that the right-hand $i^*$ is an isomorphism (since $Q$ doesn't have any rational point) to obtain 
\[
\mathrm{Sq}^2_{\omega}([h]^{n-1})=(\mathrm{Sq}^2_{\omega}\circ i^*)([h]^{n-1})=i^*([h]^n)\neq 0.\qedhere
\]
\end{proof}

We now consider the pushforward $p_*^\mathrm{I}\colon \H^n(Q,\mathrm{I}^n\otimes\omega_{Q/k})\to \mathrm{W}(k)$.

\begin{lemma}
We have $p_*^{\mathrm{I}}(\partial([h]^{n-1}))=0$
and an exact sequence
\[
0\to \Z/2\Z\xrightarrow{\partial([h]^{n-1})} \ker(p_*^{\mathrm{I}})\to \ker(p_*^{\mathrm{W}})\to 0.
\]
\end{lemma}

\begin{proof}
By construction, we have a commutative diagram
\[
\xymatrix{\H^n(Q,\mathrm{I}^n\otimes\omega_{Q/k})\ar[r]\ar[d]_-{p_*^{\mathrm{I}}} & \H^n(Q,\mathrm{I}^{n-1}\otimes\omega_{Q/k})\ar[d]^-{p_*^{\mathrm{W}}} \\
\mathrm{W}(k)\ar@{=}[r] & \mathrm{W}(k)}
\]
and $p_*^{\mathrm{I}}(\partial([h]^{n-1}))=0$ since $\partial([h]^{n-1})$ vanishes under the top horizontal arrow. Consequently, we obtain a commutative diagram
\[
\xymatrix@C=4em{0\ar[r] & \Z/2\Z\ar[r]^-{\partial([h]^{n-1})}\ar[d] & \H^n(Q,\mathrm{I}^n\otimes\omega_{Q/k})\ar[r]\ar[d]_-{p_*^{\mathrm{I}}} & \H^n(Q,\mathrm{I}^{n-1}\otimes\omega_{Q/k})\ar[d]^-{p_*^{\mathrm{W}}}\ar[r] & 0 \\
0\ar[r] & 0\ar[r] &\mathrm{W}(k)\ar@{=}[r] & \mathrm{W}(k)\ar[r] & 0}
\]
and we conclude using the snake lemma.
\end{proof}

Consider now the short exact sequence of sheaves on $Q$:
\begin{equation}\label{eqn:fundamental}
0\to 2\KM_n\to \KMW_n\otimes\omega_{Q/k}\to \mathrm{I}^n\otimes\omega_{Q/k}\to 0.
\end{equation}
The left-hand term fits into an exact sequence of sheaves 
\[
0\to {}_2\KM_n\to \KM_n\to 2\KM_n\to 0.
\]
and the $n$-fold contraction of ${}_2\KM_n$ is trivial. Thus, the morphism $\KM_n\to 2\KM_n$ induces an isomorphism $\H^n(Q,\KM_n)\to \H^n(Q,2\KM_n)$ and \eqref{eqn:fundamental} yields an exact sequence of the form
\[
\mathrm{H}^n(Q,\KM_n)\to \H^n(Q,\KMW_n\otimes\omega_{Q/k}) \to \H^n(Q,\mathrm{I}^n\otimes\omega_{Q/k})\to 0.
\]
Using the pushforward map along the projection $p\colon Q\to \Spec(k)$, we obtain a commutative diagram
\[
\xymatrix@C=1.3em{ & \mathrm{H}^n(Q,\KM_n)\ar[r]\ar[d]_-{p_*^{\mathrm{K}}} & \H^n(Q,\KMW_n\otimes\omega_{Q/k}) \ar[r]\ar[d]_-{p_*^{\mathrm{GW}}} & \H^n(Q,\mathrm{I}^n\otimes\omega_{Q/k})\ar[r]\ar[d]_-{p_*^{\mathrm{I}}} & 0\\
0\ar[r] & \Z\ar[r] & \mathrm{GW}(k)\ar[r] & \mathrm{W}(k)\ar[r] & 0}
\]
The homomorphism $\Z\to \mathrm{GW}(k)$ is the inclusion of $\Z\cdot h$ in $\mathrm{GW}(k)$, and is thus injective. On the other hand, the pushforward map $p_*^{\mathrm{K}}$ is also injective with image $2\Z$. Consequently, the left-hand top horizontal arrow is in turn injective and we obtain a morphism of short exact sequences
\begin{equation}\label{eqn:snake}
\xymatrix@C=1.3em{0\ar[r] & \mathrm{H}^n(Q,\KM_n)\ar[r]\ar[d]_-{p_*^{\mathrm{K}}} & \H^n(Q,\KMW_n\otimes\omega_{Q/k}) \ar[r]\ar[d]_-{p_*^{\mathrm{GW}}} & \H^n(Q,\mathrm{I}^n\otimes\omega_{Q/k})\ar[r]\ar[d]_-{p_*^{\mathrm{I}}} & 0\\
0\ar[r] & \Z\ar[r] & \mathrm{GW}(k)\ar[r] & \mathrm{W}(k)\ar[r] & 0.}
\end{equation}
The snake lemma then yields an exact sequence
\[
0\to \ker(p_*^{\mathrm{GW}})\to \ker(p_*^{\mathrm{I}})\to \Z/2\Z
\]
\begin{lemma}
The image of $\partial([h]^{n-1})$ under the homomorphism 
\[
\ker(p_*^{\mathrm{I}})\to \Z/2\Z
\]
is nontrivial.
\end{lemma}

\begin{proof}
Let $Q_1=Q\cap \{x_4=\ldots=x_{n+2}=0\}=Q\cap \mathbb{P}^2\subset \mathbb{P}^{n+1}$ be a smooth conic section. The inclusion $j\colon Q_1\to Q$ gives a commutative diagram
\[
\xymatrix{\H^0(Q_1,\overline{\mathrm{W}})\ar[r]^-{\partial'}\ar[d]_-{j_*} & \H^1(Q_1,\mathrm{I}\otimes\omega_{Q_1/k})\ar[d]^-{j_*} \\
\H^{n-1}(Q,\overline{\mathrm{I}}{}^{n-1})\ar[r]_-\partial & \H^n(Q,\mathrm{I}^n\otimes\omega_{Q/k}).}
\]
We have $[h]^{n-1}=j_*(\overline 1)$ and consequently $\partial([h]^{n-1})=j_*(\partial'(\overline 1))$. We know from \cite[Lemma 3.1]{Fas13} that $\partial'(\overline 1)=e(\omega_{Q_1/k}^\vee)=e(\mathcal{O}(1))\in  \H^1(Q_1,\mathrm{I}\otimes\omega_{Q_1/k})$. If we denote by $\pi\colon Q_1\to \Spec(k)$ the projection, we obtain that $0=p_*^{\mathrm{I}}(\partial([h]^{n-1}))=\pi_*(e(\mathcal{O}(1)))$. On the other hand, there is an obvious lift of $e(\mathcal{O}(1))$ in $\H^1(Q_1,\KMW_1\otimes\omega_{Q_1/k})$, namely the (Milnor--Witt theoretic) Euler class that we still denote by $e(\mathcal{O}(1))$. Its pushforward under $\pi_*^{\mathrm{GW}}\colon \H^1(Q_1,\KMW_1\otimes\omega_{Q_1/k})\to \mathrm{GW}(k)$ has trivial Witt class (since $\pi_*(e(\mathcal{O}(1)))=0$) and rank $2$ (since its rank is the pushforward of $c_1(\mathcal{O}(1))$, yielding $\pi_*^{\mathrm{GW}}( e(\mathcal{O}(1)))=h$.
\end{proof}

\begin{proof}[Proof of Theorem \ref{thm:pushforwardequivalence}]
Consider the commutative diagram
\[
\xymatrix{ & \Z/2\Z\ar[d]_-{\partial([h]^{n-1})}\ar@{-->}[rd] & \\ 
\ker(p_*^{\mathrm{GW}})\ar[r]\ar@{-->}[rd] & \ker(p_*^{\mathrm{I}})\ar[r]\ar[d] & \Z/2\Z \\
 & \ker(p_*^{\mathrm{W}}) & }
\]
in which the row and column are short exact and the top diagonal map is an isomorphism. A straightforward diagram chase shows that the bottom diagonal arrow is also an isomorphism.
\end{proof}

Using Theorem \ref{thm:explicitkernel}, we obtain the following corollary.

\begin{corollary}\label{cor:explicitkernel}
Suppose that $1\leq n\leq 7$. If $n\neq 4$, there is an exact sequence
\[
\mathrm{W}^{n-5}(Q,\omega_{Q/k})\xrightarrow{e}\H^{n-5}(Q,\mathrm{W}\otimes \omega_{Q/k})\xrightarrow{\partial} \CHW_0(Q)\xrightarrow{p_*}\mathrm{GW}(k)
\]
where $\partial$ is induced by the differential at page $5$ of the Gersten--Witt spectral sequence and $e$ is an edge map. If $n=4$, then there is an exact sequence of the form
\[
0\to \ker(sp^*)/\ker(p^*)\to\CHW_0(Q)\xrightarrow{p_*}\mathrm{GW}(k).
\]
\end{corollary}

We now prove that there are some examples of quadrics of dimension $4$ where $ \ker(sp^*)/\ker(p^*)$ is nontrivial. We are indebted to Stephen Scully for pointing out the case of Albert forms.

\begin{lemma}\label{lem:Albert}
Suppose that $q$ is an anisotropic Albert form, i.e. $q=\langle -a,-b,ab,c,d,-cd\rangle$ for some $a,b,c,d\in k^\times$ and the biquaternion algebra $A:=(a,b)_k\otimes_k(c,d)_k$ is a division algebra (\cite[IV,Theorem 16.5]{KMRT98}). Then the map 
\[
p^*\colon \mathrm{W}(k)\to \mathrm{W}^0(Q,\omega_{Q/k})
\]
is injective.
\end{lemma}

\begin{proof}
Using the exact sequence of Corollary \ref{cor:exactsequences}, it suffices to prove that $\mathrm{W}^0(C_0(q)_\sigma)=0$. We know from \cite[IV, proof of Proposition 16.3]{KMRT98} that 
\[
(C_0(q),\sigma)\simeq (A\times A^{op},\epsilon)
\]
as algebras with involutions, where $\epsilon$ is the exchange involution. This immediately implies that $\mathrm{W}^0(C_0(q)_\sigma)=0$ \cite[Example~2.4]{Fir23}.
\end{proof}

As a consequence, we see that the exact sequence of Corollary \ref{cor:explicitkernel} becomes
\[
0\to \ker(sp^*)\to\CHW_0(Q)\xrightarrow{p_*}\mathrm{GW}(k).
\]

\begin{proposition}\label{prop:counterexample}
If $k=\mathbb{C}(T_1,T_2,T_3,T_4)$, $q:=\langle T_1,T_2,T_1T_2,T_3,T_4,T_3T_4\rangle$ and $Q$ is the associated quadric, the pushforward homomorphism $p_*\colon \CHW_0(Q)\to \mathrm{GW}(k)$ is not injective.
\end{proposition}

\begin{proof}
Since $-1$ is a square in $k$, $q$ is isometric to $\langle -T_1,-T_2,T_1T_2,-T_3,-T_4,T_3T_4\rangle$, which is a subform of the 4-fold Pfister form $\langle\langle T_1,T_2,T_3,T_4\rangle\rangle$. The latter is anisotropic, as a consecutive use of the $T_i$-adic valuation easily shows. It follows that $q$ is anisotropic and the previous lemma shows that $p^*\colon \mathrm{W}(k)\to \mathrm{W}(Q)$ is injective. It remains to observe that the map $sp^*\colon \mathrm{W}(k)\to \mathrm{W}(k(Q))$ is not injective. Now, $q$ and therefore $\langle\langle T_1,T_2,T_3,T_4\rangle\rangle$ become isotropic over $k(Q)$, showing that the latter is hyperbolic. Thus $sp^*(\langle\langle T_1,T_2,T_3,T_4\rangle\rangle)=0$, but its class in $\mathrm{W}(k)$ is nontrivial. 
\end{proof}

\section{Suslin matrices and Euler classes}\label{sec:SuslinEuler}

In this section we give an example of an affine quadric of dimension $4$ for which the converse of \cite[Theorem~4.1.2]{AL25} does not hold. We use the Suslin matrices to study the Euler class of the tangent bundle and apply the results of the previous section to show that it does not vanish.

\begin{definition}
    Following \cite[\S~5]{Sus77}, for $x:=(x_1,x_2,\hdots,x_n)$, $y:=(y_1,y_2,\hdots, y_n)$,
    we define matrices $\alpha_n(x,y)$  of size $2^{n-1}$ inductively as
    \begin{gather*}
    \alpha_1(x_1,y_1):=x_1,\quad x':=(x_2,x_3,\hdots,x_n),\quad y':=(y_2,y_3,\dots,y_n),\\
    \alpha_n(x,y)=\begin{pmatrix}
        x_1\cdot \Id_{2^{n-2}} & \alpha_{n-1}(x',y') \\
        -\alpha_{n-1}(y',x')^T & y_1\cdot \Id_{2^{n-2}} 
    \end{pmatrix}.
    \end{gather*}
    Continuing to follow \cite[\S~5]{Sus77}, define $I_1:=1$ and
    \[
    I_n:=\begin{cases}
        \begin{pmatrix}
            0 & I_{n-1} \\
            -I_{n-1} & 0 
        \end{pmatrix}
        & \text{if } n=2m,\\
        \begin{pmatrix}
            I_{n-1} & 0 \\
            0 & -I_{n-1} 
        \end{pmatrix}
        & \text{if } n=2m+1.
    \end{cases}
    \]
    These matrices satisfy $I_n^T=I_n^{-1}=(-1)^{\frac{n(n-1)}{2}}\cdot I_n$ and $\det I_n=1$ by \cite[Lemma~5.2]{Sus77}.
    
    Let $a:=(a_1,a_2,\hdots, a_n)$, we denote $\diag(a):=\diag(a_1,a_2,\hdots,a_n)$ and 
    \[
    \Sus{a}:=\alpha_n(x,x\cdot \diag(a))\cdot I_n.
    \]
    Recall that \cite[Lemma~5.3]{Sus77} yields
    \[
    \Sus{a}^T = \begin{cases}
        \Sus{a} & \text{if $n=4m+1$},\\
        -\Sus{a} & \text{if $n=4m+3$}.
    \end{cases}
    \]
    Finally, define diagonal matrices $\Lev{a}$ and $\Lodd{a}$ of size $2^{n-1}$ inductively as 
    \begin{gather*}
    \Lev{a_n}:=1,\quad \Lodd{a_n}:=a_n, \quad a':=(a_2,a_3,\hdots,a_n),\\
    \Lev{a}=\begin{pmatrix}
        \Lev{a'} & 0 \\
        0 & a_1\Lodd{a'}
    \end{pmatrix},\quad \Lodd{a}=\begin{pmatrix}
        a_1\Lev{a'} & 0 \\
        0 & \Lodd{a'}
    \end{pmatrix}.
    \end{gather*}
    The matrix $\Lev{a}$ (resp. $\Lodd{a}$) is diagonal with the entries being products of even (resp. odd) number of $a_i$-s ordered in a particular way.
\end{definition}

\begin{lemma}\label{lem:suslin_preliminary}
    Let $a:=(a_1,a_2,\hdots, a_n)$ for $a_i\in k^\times$. The following holds:
    \begin{enumerate}
        \item  For $\Delta:=\prod_{i=1}^na_i$ one has
    \begin{gather*}
    I_n \cdot \Lev{a} = \begin{cases}
        \Delta \cdot (\Lev{a})^{-1}  \cdot I_n, & n=2m,\\
        \Delta \cdot (\Lodd{a})^{-1}  \cdot I_n, & n=2m+1,
    \end{cases}\\
        I_n \cdot \Lodd{a} = \begin{cases}
        \Delta \cdot (\Lodd{a})^{-1}  \cdot I_n, & n=2m,\\
        \Delta \cdot (\Lev{a})^{-1} \cdot  I_n, & n=2m+1.
    \end{cases}
    \end{gather*}
    \item Let $x=(x_1,x_2,\hdots,x_n)$. Then
        \[
        \alpha_n(x\cdot \diag(a),x)=\Lev{a}^{-1}\cdot \alpha_{n}(x,x\cdot \diag(a)) \cdot \Lodd{a}.
        \]
    \end{enumerate}
\end{lemma}
\begin{proof}
    Follows immediately from the inductive definition.
\end{proof}

\begin{lemma}\label{lem:Suslin_identity}
    Let $a:=(a_1,a_2,\hdots, a_n)$, denote
    $\Delta:=\prod_{i=1}^{n} a_i$ and $q:=\sum_{i=1}^{n}a_ix_i^2$.
    Then
    \[
    \Delta \cdot \Sus{a}^T (\Lev{a})^{-1} \Sus{a} = \begin{cases}
            q \cdot \Lodd{a}, & n=2m,\\
        q \cdot \Lev{a}, & n=2m+1,
    \end{cases}
    \]
\end{lemma}
\begin{proof}
    Let $x:=(x_1,x_2,\hdots,x_n)$ and $N:=2^{n-1}$, then \cite[Lemma~5.1]{Sus77} yields
    \[
    \alpha_n(x\cdot \diag(a),x)^T\cdot \alpha_n(x,x\cdot \diag(a)) =q\cdot \Id_{N}.
    \]
    Applying Lemma~\ref{lem:suslin_preliminary}.2 we obtain
    \[
    \alpha_n(x,x\cdot \diag(a))^T \cdot \Lev{a}^{-1} \cdot \alpha_n(x,x\cdot \diag(a)) = q\cdot \Lodd{a}^{-1}.
    \]
    Multiplying both sides with $I_n^T$ on the left and $I_n$ on the right we obtain
    \[
    \Sus{a}^T \cdot \Lev{a}^{-1} \cdot \Sus{a} = q\cdot I_n^T\cdot\Lodd{a}^{-1}\cdot I_n.
    \]
    The claim follows from Lemma~\ref{lem:suslin_preliminary}.1 and the identity $I_n^TI_n=\id_N$. 
\end{proof}

\begin{corollary} \label{cor:Suslin_constant}
    Let $a:=(a_1,a_2,\hdots, a_{4m+1})$ with $a_i\in k^\times$ and denote
    \[
    \Delta:=\prod_{i=1}^{4m+1} a_i,\quad q:=\sum_{i=1}^{4m+1} a_ix_i^2,\quad R:=k[x_1,x_2,\hdots,x_{4m+1}]/(q-1).
    \]
    Suppose that $\Delta=-d^2$ for some $d\in k^\times$. Then
    \[
    \langle 2d\cdot \Sus{a}\rangle=\sum_{j=0}^{2m}\Lambda^{2j}\langle a_1,a_2,\hdots,a_{4m+1}\rangle
    \]
    in $\GW(R)$, where $\langle 2d\cdot \Sus{a}\rangle$ is the quadratic form given by the symmetric matrix $2d\cdot \Sus{a}$.
\end{corollary}
\begin{proof}
    Let $N:=2^{4m}$ and $  A:=\Id_{N}+d\cdot (\Lev{a})^{-1}\Sus{a}$.
    By Lemma~\ref{lem:Suslin_identity}, we have 
    \[
    A\cdot (\Id_{N}-d\cdot (\Lev{a})^{-1}\Sus{a}) = (1+q)\cdot\Id_{N},
    \]
 thus $A$ is invertible over $R$.
    Since $\Sus{a}$ is symmetric, $A^T=\Id_{N}+d\cdot \Sus{a}(\Lev{a})^{-1}$. Lemma~\ref{lem:Suslin_identity} furthermore yields
    \begin{multline*}
    A^T \Lev{a} A = (\Id_{N}+d\cdot \Sus{a}(\Lev{a})^{-1})\cdot \Lev{a} \cdot (\Id_{N}+d\cdot(\Lev{a})^{-1}\Sus{a}) = \\
    2d\cdot \Sus{a} + (1-q)\cdot \Lev{a}.
    \end{multline*}
    Thus $2d\cdot \Sus{a}$ is congruent to $\Lev{a}$ over $R$, the claim follows.
\end{proof}

\begin{remark}
    Corollary~\ref{cor:Suslin_constant} admits the following Clifford algebra interpretation. Let $C(q)$ be the Clifford algebra equipped with the standard symmetric bilinear form $\langle u,v\rangle:={\frac{1}{2^n}}\tr (\sigma(u)v)$, where $\sigma$ is the reversal antiautomorphism and ${\frac{1}{2^n}}\tr (w)$ is the scalar part of $w\in C(q)$. Put 
    \[
    C_R:=R\otimes_k C(q),\quad x:=x_1e_1+x_2e_2+\hdots +x_ne_n,
    \]
    for the standard basis $e_1,\hdots,e_n$ of $k^n$, 
    and extend $\langle-,-\rangle$ to $C_R$. Consider the homomorphism 
    \[
    J\colon C_R\to C_R,\quad v\mapsto x v.
    \]
    
    Put
    \[
    e:=e_1e_2\hdots e_n,\quad \epsilon:=de^{-1},
    \]
    where $\Delta=-d^2$ as above.
    Since $n=4m+1$, $\epsilon$ is central in $C(q)$ as well as in $C_R$, $\sigma(\epsilon)=\epsilon$, and $\epsilon^2=-1$. Then $(\epsilon J)^2=-\id$ and $\langle \epsilon Ju,v\rangle =\langle u,\epsilon Jv\rangle$, hence
    \begin{equation} \label{eq:Suslin_const_Cliffors}
    \langle(1+\epsilon J)u,(1+\epsilon J)v\rangle=\langle u,v\rangle+2\langle u,\epsilon Jv\rangle+\langle u,(\epsilon J)^2v\rangle =2 \langle u,\epsilon Jv\rangle. 
    \end{equation}
    Furthermore, $1+\epsilon J$ is invertible with the inverse given by $\frac{1}{2}(1-\epsilon J)$.

    Let $\mc{B}$ be the canonical basis of the even part $C_0(q)$ of $C(q)$ in which the Gram matrix of $\langle -,-\rangle$ restricted $C_0(q)$ is given by $\Lev{a}$. 
    One may show by induction (cf. \cite[\S~3]{Sye26}) that in the basis $\mc{B}$ the Suslin matrix $\Sus{a}$ satisfies
    \[
    [e^{-1}J]_\mc{B}=(\Lev{a})^{-1} \Sus{a}.
    \]
    Taking the Gram matrices in the equation~\eqref{eq:Suslin_const_Cliffors} we obtain
    \[
    ([1+\epsilon J]_\mc{B})^T\cdot \Lev{a}\cdot [1+\epsilon J]_\mc{B}=2d\cdot \Sus{a}.
    \]
\end{remark}

Till the end of the section we put
\begin{gather*}
k:=\C(T_1,T_2,T_3,T_4),\quad q:=T_1T_2x_1^2+T_2x_2^2+T_1T_3x_3^2+T_1T_4x_4^2+T_1T_3T_4x_5^2,\\ \tilde{q}:=T_1x_0^2+T_2x_1^2+T_1T_2x_2^2+T_3x_3^2+T_4x_4^2+T_3T_4x_5^2, \\
Q:=\{\tilde{q}=0\}\subseteq \PP^5,\quad Q^o:=Q\setminus\{x_0=0\},\quad Q^o\cong \Spec R, \quad R:=k[x_1,\ldots,x_5]/(q-1).
\end{gather*}

\begin{corollary}\label{cor:explicitcomputation}
$\langle S(T_1T_2,T_2,T_1T_3,T_1T_4,T_1T_3T_4)\rangle=\langle\langle T_1,T_2,T_3,T_4\rangle\rangle$ in $\mathrm{GW}(R)$.
\end{corollary}

\begin{proof}
We have $\Delta(q)=T_1^4T_2^2T_3^2T_4^2$, which satisfies the condition of Corollary~\ref{cor:Suslin_constant} with $d=T_2T_3T_4$ up to the square factor $(\sqrt[4]{-1}\cdot T_1)^2$. Now, a direct computation using $\sqrt{-1}\in k$ shows that 
\[
\langle\Lev{T_1T_2,T_2,T_1T_3,T_1T_4,T_1T_3T_4}\rangle=\langle\langle T_1,T_2,T_3,T_4\rangle\rangle
\]
and Corollary~\ref{cor:Suslin_constant} yields 
\[
\langle S(T_1T_2,T_2,T_1T_3,T_1T_4,T_1T_3T_4)\rangle=\langle T_2T_3T_4\rangle\langle\langle T_1,T_2,T_3,T_4\rangle\rangle.
\]
To conclude, we observe that $\mathrm{I}^5(k)=0$, and thus 
\[
\langle T_2T_3T_4\rangle\langle\langle T_1,T_2,T_3,T_4\rangle\rangle=\langle\langle T_1,T_2,T_3,T_4\rangle\rangle. \qedhere
\]
\end{proof}

\begin{corollary}\label{cor:eulerclassopenquadric}
The image of the Euler class $e(T_Q)$ under the composite
\[
\CHW_0(Q)\to \CHW^4(Q)\to \mathrm{GW}(Q)\to  \mathrm{GW}(Q^o), 
\]
given respectively by the canonical choice of the square root of $\omega_{Q/k}\cong \struct(-4)\cong \struct(-2)^{\otimes 2}$, the edge map in the Grothendieck--Witt spectral sequence and the localization homomorphism to $Q^o$, is $\langle S(T_1T_2,T_2,T_1T_3,T_1T_4,T_1T_3T_4)\rangle-8h$ (up to a unit in $\mathrm{W}(k)$). 
\end{corollary}

\begin{proof}
Localization induces a commutative diagram
\[
\xymatrix{\CHW_0(Q)\ar[r]\ar[d] & \CHW^4(Q)\ar[r]\ar[d] & \mathrm{GW}(Q)\ar[d] \\
\CHW_0(Q^o)\ar[r] & \CHW^4(Q^o)\ar[r] & \mathrm{GW}(Q^o)}
\]
and it suffices thus to prove that the image of $e(T_{Q^o})$ under the edge map $\CHW^4(Q^o)\to \mathrm{GW}(Q^o)$ is $\langle\langle T_1,T_2,T_3,T_4\rangle\rangle -8h$ (up to a unit). By \eqref{eqn:etaphi}, the Euler class is given by the cycle $\langle {x_1} \rangle\otimes  ({dx}_2\wedge {dx}_3\wedge\hdots\wedge {dx}_{5})^\vee$, and we consider the closed subscheme $Z:=\{x_2=\ldots=x_5=0\}\subset Q^o$. The global section $x_1$ is invertible on $Z$ and the Euler class is the image of $\langle x_1\rangle\in \CHW^0(Z)$ under the pushforward $\CHW^0(Z)\to \CHW^4(Q^o)$ (obtained via the trivialization of the normal bundle given by $x_2,\ldots,x_5$). We have a commutative diagram
\[
\xymatrix{\CHW^0(Z)\ar[r]\ar[d] & \mathrm{GW}^0(Z)\ar[d]  \\
\CHW^4(Q^o)\ar[r] & \mathrm{GW}^0(Q^o) }
\]
where the horizontal maps are edge maps in the corresponding spectral sequences, and the vertical maps are pushforwards along the embedding $Z\to Q^o$. Using this, we see that the image of the Euler class under the edge map is the pushforward of $\langle x_1\rangle\in \mathrm{GW}^0(Z)$.

This computation can be done universally: Let $Q_9\subset\A^{10}_k$ be the smooth affine quadric with coordinate ring 
\[
k[Q_9]=k[u_1,\ldots,u_5,v_1,\ldots,v_5]/(\sum_{i=1}^5 u_iv_i-1),
\] 
and let $U$ be the universal bundle of rank $4$ on $Q_9$, i.e. the kernel of the unimodular row given by the five last coordinates.

The choice of the unimodular rows $(x_1,\ldots,x_5)$ and $(T_1T_2x_1,T_2x_2,T_1T_3x_3,T_1T_4x_4,T_1T_3T_4x_5)$ yields a morphism $\varphi\colon Q^o\to Q_9$ such that $\varphi^*(U)=T_{Q^o}$. If $T:=\{u_2=\ldots=u_5=0\}\subset Q_9$, we have $T\simeq \mathbb{G}_m\times \A^4$ and a Cartesian square
\[
\xymatrix{Z\ar[r]\ar[d]_-{\varphi_{\vert Z}} & Q^o \ar[d]^-\varphi \\
T\ar[r] & Q_9}
\]
in which the horizontal maps are regular closed immersions of codimension $4$ satisfying $\varphi_{\vert Z}^*N_{Q_9}T=N_{Q^o}Z$. Since $\varphi_{\vert Z}^*(\langle u_1\rangle)=\langle x_1\rangle\in \mathrm{GW}^0(Z)$, we may use the excess intersection formula \cite[Theorem 32]{Fas09} to reduce the problem to the case of $Q_9$, i.e. the computation of the image of the Euler class $e(U)$ under the edge homomorphism
\[
\CHW^4(Q_9)\to \mathrm{GW}^0(Q_9).
\]

By \cite[\S 3.3]{Fas11b}, the free $\mathrm{W}(k)$-module $\CHW^4(Q_9)$ admits $\eta[u_1]\otimes  ({du}_2\wedge {du}_3\wedge\hdots\wedge {du}_{5})^\vee $ as generator, and Lemma \ref{lem:etaEulerclass} shows that this is the Euler class of $U$. On the other hand, the reduced group $\widetilde{\mathrm{GW}}^0(Q_9)$ (i.e. the kernel of the pullback map $\mathrm{GW}^0(Q_9)\to \mathrm{GW}(k)$ given by a rational point) is shown in \cite[Proposition 3.4.3]{AF17} to be the free $\mathrm{W}(k)$-module generated by the class of $[\alpha_n(u,v)\cdot I_n]-8h$, where $u:=(u_1,\ldots,u_5), v:=(v_1,\ldots,v_5)$.  Using \cite[Proposition 1.1.5]{AF22}, we deduce that the edge map $\CHW^4(Q_9)\to \widetilde{\mathrm{GW}}^0(Q_9)$ is an isomorphism, and thus maps $e(U)$ to $[\alpha_n(u,v)\cdot I_n]-8h$ up to a unit in $\mathrm{W}(k)$. We conclude using the fact that $\varphi^*([\alpha_n(u,v)\cdot I_n]-8h)=[S(a)]-8h$.
\end{proof}

\begin{remark}
In fact, we could consider the quadric $Q_9$ as defined over $\Z$, and the above arguments show that the image of the Euler class is equal to the Suslin matrix up to a hyperbolic summand and multiplication by a unit in $\mathrm{W}(\Z)$, i.e. by $\langle \pm 1\rangle$. 
\end{remark}

We can finally prove the main result of the section. 

\begin{theorem}
The tangent bundle to $Q^o$ has nontrivial Euler class.
\end{theorem}

\begin{proof}
In view of Corollary \ref{cor:eulerclassopenquadric} the image of the Euler class $e(T_{Q^o})$ under the edge homomorphism is (up to a unit) $\langle S(T_1T_2,T_2,T_1T_3,T_1T_4,T_1T_3T_4)\rangle-8h$. We know from Corollary \ref{cor:explicitcomputation} that this class is just the class of $\langle\langle T_1,T_2,T_3,T_4\rangle\rangle -8h$ and it suffices to show that its class in $\mathrm{W}(Q^o)$ is nontrivial to conclude. However, by the proof of Proposition~\ref{prop:counterexample}, this class, seen as an element of $\mathrm{W}(Q)$, is nontrivial and in the kernel of $p_*\colon \mathrm{W}(Q)\to \mathrm{W}(k)$. Now, we consider the diagram
\[
\xymatrix{ 
 & \mathrm{W}^3(Q_{\infty},\omega_{Q_{\infty}/k})\ar[d]^-{j_*}\ar[rd]^-{p_*j_*} & \\
\ker(p_*)\ar[r]^-{\iota}\ar[rd]_-{s\iota} &  \mathrm{W}(Q)\ar[r]^-{p_*}\ar[d]^-s &\mathrm{W}(k) \\
 & \mathrm{W}(Q^o) & }
\]
in which $j,s$ are the homomorphisms in the long exact sequence of localization and $\iota$ is the inclusion. As $p_*j_*$ is injective by Corollary~\ref{cor:exactsequences}, a diagram chase shows that $s\iota$ as well. Thus the class of $\langle\langle T_1,T_2,T_3,T_4\rangle\rangle$ in $\mathrm{W}(Q^o)$ is nontrivial.
\end{proof}

\appendix

\section{On $\K$-cohomology groups of quadrics}\label{sec:Rosttheorem}

In this section we show that for a projective quadric $Q$ of  dimension $n$ the groups $\H^{n}(Q,\K^\M_{n+1})$, $\H^{n-1}(Q,\kM_{n})$ and $\Ch^{n-1}(Q)$ are generated by cycles of particular kind. Theorem~\ref{thm:rost} below is a particular case of a more general result of Rost (see \cite[Remark~4.6]{Mer00}). Since the general result remains unpublished we provide a proof for the case of quadrics that we need in the paper. The proof is based on a sketch of an argument suggested by an LLM to the first author when he was trying to locate a published version of the aforementioned result of Rost. 

\begin{lemma} \label{lem:deg2_generation}
    Let $Q$ be a smooth projective quadric of dimension $n$ over $k$. Then $\H^n(Q,\K^\M_{n+1})$ is generated by cycles $\{u\}_z$, $z\in Q_{(0)}$, $u\in k(z)^\times$ with $[k(z):k]\le 2$.
\end{lemma}
\begin{proof}
    If $n=0$ the claim is trivial, so suppose $n\ge 1$. If $Q(k)\neq\emptyset$ then the pushforward $\pi_*\colon \H^n(Q,\K^\M_{n+1}) \to \K_1^\M(k)$ is an isomorphism by Corollary~\ref{cor:stable_birational_invariant}, so $\H^n(Q,\K^\M_{n+1})$ is generated by $\{u\}_z$, $u\in k^\times$, for a rational point $z\in Q$. If $Q(k)=\emptyset$ then the claim follows from \cite[Theorem~6.2, Remark~6.4]{CM01} and the definition of $\hat{\omega}_\varphi$ from \cite[Theorem~12]{Ros88}.
\end{proof}

\begin{lemma} \label{lem:1-cycles_on_isotropic_quadric}
    Let $Q$ be a smooth projective quadric of dimension $n\ge 2$ over $k$ and suppose $Q(k)\neq \emptyset$. Then 
    \[
    \H^{n-1}(Q,\K^\M_{n})=\sum_{\substack{\text{$L/k$ finite}\\ i(q_L)\ge 2}} (\pi_{L/k})_*\left( L^\times\cdot [\ell_L]\right),
    \]
    where $\pi_{L/k}$ is the projection $Q_L\to Q$ and $[\ell_L]\in \CH^{n-1}(Q_L)$ is the class of a line.
\end{lemma}
\begin{proof}
    The proof is based on the well-known geometric construction used to decompose the motive of an isotropic quadric, see e.g. \cite[Example~66.7]{EKM08}.

    Pick a rational point $x_0\in Q(k)$ and denote by $H$ the hyperplane tangent to $Q$ at $x_0$. Consider the localization sequence
    \[
    \H^{n-2}((Q\cap H)\setminus\{x_0\},\K^\M_{n-1}) \xrightarrow{i_Z} \H^{n-1}(Q\setminus \{x_0\},\K^\M_{n}) \to \H^{n-1}(Q\setminus H,\K^\M_{n}),
    \]
    we have $Q\setminus H\cong \A^n$ by the stereographic projection, hence $i_Z$ is surjective. Furthermore, the same projection gives rise to an affine bundle $p\colon(Q\cap H)\setminus\{x_0\} \to Q'$, where $Q'$ is the projective quadric corresponding to $q'$ such that $q=\langle 1,-1\rangle+q'$. It follows that $\H^{n-1}(Q\setminus \{x_0\},\K^\M_{n})$ is generated by the classes $(i_x)_*(\{u\})$, $u\in k(x)^\times$, $x\in Q'_{(0)}$, where $i_x\colon \A^1_{k(x)}\to Q\setminus \{x_0\}$ is the embedding of the fiber $p^{-1}(x)$. By localization sequence the pullback 
    \[
    j^*\colon \H^{n-1}(Q,\K^\M_{n}) \to\H^{n-1}(Q\setminus \{x_0\},\K^\M_{n})
    \]
    is an isomorphism. For every $x\in Q'_{(0)}$ we have a cartesian square
    \[
    \xymatrix{
    \A^1_{k(x)} \ar[r] \ar[d]_{i_x} & \PP^1_{k(x)} \ar[d]^{\bar{\iota}_x} \\
    Q\setminus \{x_0\} \ar[r]^(0.6)j& Q
    }
    \]
    By base change for $u\in k(x)^\times$ we have $(i_x)_*(\{u\})=j^*(\bar{\iota}_x)_*(\{u\})$, hence the classes $(\bar{\iota}_x)_*(\{u\})$ generate $\H^{n-1}(Q,\K^\M_{n})$. Since $q'_{k(x)}$ is isotropic, $i(q_{k(x)})\ge 2$. Finally, $\bar{\iota}_x$ factors as the composition $\PP^1_{k(x)}\subseteq Q_{k(x)} \to Q$. The claim follows.
\end{proof}

The next lemma seems to be folklore,  cf. \cite[Proposition~2.2]{BR20} for a split quadric surface and the MathOverflow answer \cite{MO13} for a split dimension $3$ quadric. The proof below is basically a straightforward extension of the MathOverflow answer.
\begin{lemma}\label{lem:Hilb_Q}
    Let $q$ be a regular quadratic form on a vector space $V$ over $k$, and let $Q\subseteq \PP(V)$ and $\mr{OGr(2,q)}\subseteq\Gr(2,V)$ be the respective projective quadric and the Grassmannian of totally isotropic planes. Denote by 
    \[
    Q^{[2]}:= (\mr{Bl}_{\Delta} (Q\times Q))/(\Z/2)
    \]
    the Hilbert scheme of length $2$ subschemes. Then there is an isomorphism
    \[
    Q^{[2]} \cong \mr{Bl}_{\mr{OGr(2,q)}} \Gr(2,V).
    \]
\end{lemma}
\begin{proof}
Let $\mc{T}$ be the tautological rank $2$ bundle over $\Gr(2,V)$. The quadratic form $q\in \Sym^2 V^\vee$ gives rise to a section $s_q$ of $\Sym^2\mc{T}^\vee$ whose zero locus is $\mr{OGr}(2,q)$. Recall (see e.g. \cite[\S4.1]{Goe94}) that
\[
\PP(V)^{[2]}\cong \PP(\Sym^2\mc{T}^\vee).
\]
A point of $\PP(\Sym^2\mc{T}^\vee)$ corresponds to a line
$\ell=\PP(W)\subseteq\PP(V)$ together with a degree two effective divisor $D\subseteq\ell$ given by $f\in \PP(\Sym^2W^\vee)$. Under this identification, $Q^{[2]}\subseteq \PP(V)^{[2]}\cong \PP(\Sym^2\mc{T}^\vee)$
consists of pairs $(\ell,D)$ such that $D\subseteq Q$, which is equivalent to $q|_\ell\in f$. Thus $Q^{[2]}$ is the zero scheme of the section of $p^*(\Sym^2\mc{T}^\vee)/\struct_{\PP(\Sym^2\mc{T}^\vee)}(-1)$
induced by $s_q$, here $p\colon \PP(\Sym^2\mc{T}^\vee)\to \Gr(2,V)$ is the projection.

If $\dim Q\leq 1$ then the zero locus of $s_q$ is empty, hence $Q^{[2]}\cong \Gr(2,V)$. Otherwise, the subvariety $\mr{OGr}(2,q)\subseteq \Gr(2,V)$ is smooth of codimension three and the three local components $s_1,s_2,s_3$ of $s_q$ form a regular sequence. By the above, $Q^{[2]}$ is locally given in $\Gr(2,V)\times\mathbb P^2$ by the relation
\[
(s_1,s_2,s_3)\in f,
\]
or equivalently by the equations
\[
s_if_j-s_jf_i=0,\quad 1\leq i<j\leq 3, \qquad f=[f_1:f_2:f_3].
\]
Since $(s_1,s_2,s_3)$ is a regular sequence, these equations give the standard description of the blow-up along the ideal $(s_1,s_2,s_3)$.
\end{proof}

The following is a particular case of an unpublished result by Rost, see \cite[Remark~4.6]{Mer00}.
\begin{theorem}[{Rost}] \label{thm:rost}
    Let $Q$ be a smooth projective quadric of dimension $n$ over $k$. Then
    \[
    \H^{n-1}(Q,\kM_n) = \sum_{\text{$L/k$ finite}} (\pi_{L/k})_* \left(  \kM_1(L) \cdot \H^{n-1}(Q_L,\kM_{n-1})\right)
    \]
    where $\pi_{L/k}$ is the projection $Q_L\to Q$.
\end{theorem}
\begin{proof}
    Recall that $\H^{n-1}(Q,\KM_n)\to \H^{n-1}(Q,\kM_n)$ is surjective since $\CH_0(Q)$ is torsion-free, hence the case of $n=1$ is clear, and if $Q(k)\neq\emptyset $ then the claim follows from Lemma~\ref{lem:1-cycles_on_isotropic_quadric}. So we assume that $Q(k)=\emptyset $.

    Let $
    Q^{[2]}:= (\mr{Bl}_{\Delta} (Q\times Q))/(\Z/2)$
    be the Hilbert scheme of length $2$ subschemes of $Q$. In the notation of Lemma~\ref{lem:Hilb_Q} we have
    \[
    Q^{[2]}\cong \mr{Bl}_{\mr{OGr}(2,q)} \Gr(2,V).
    \]
    Denote by $\mc{T}$ the tautological bundle over $\Gr(2,V)$.
    Since $\mr{OGr}(2,q)$ is the regular zero locus of the section $\Sym^2\mc{T}^\vee$ given by the quadratic form defining $Q$, we have $\mc{N}\cong \Sym^2\mc{T}^\vee|_{\mr{OGr}(2,q)}$ for the normal bundle $\mc{N}$ of the embedding $\mr{OGr}(2,q)\to\Gr(2,V)$. Blow-up formula yields 
    \begin{multline*}
    \H^{2n-1}(Q^{[2]},\kM_{2n})=r^*\H^{2n-1}(\Gr(2,V),\kM_{2n})\oplus s_*h^*(\H^{2n-2}(\mr{OGr}(2,q),\kM_{2n-1}))\oplus\\
    \oplus s_*(h^*\H^{2n-3}(\mr{OGr}(2,q),\kM_{2n-2})\cdot c_1),
    \end{multline*}
    where $r\colon  Q^{[2]}\to\Gr(2,V)$ is the projection, $s\colon \PP(\mc{N})\to Q^{[2]}$ is the embedding of the exceptional divisor of the blow-up, $h\colon\PP(\mc{N})\to \mr{OGr}(2,q)$ is the projection and $c_1=c_1(\struct_{\PP(\mc{N})}(1))$. The middle summand vanishes by dimension reasons.
    For the first summand cellular decomposition yields that
    \[
    \H^{2n-1}(\Gr(2,V),\kM_{2n})=\kM_1(k)\cdot \H^{2n-1}(\Gr(2,V),\kM_{2n-1}).
    \]
    Since $\dim \mr{OGr}(2,q)=2n-3$, the last summand is generated by the elements $i_*(\{u\}\cdot [\ell])$ where $u\in k(x)^\times$, $x\in \mr{OGr}(2,q)$ is a closed point, $[\ell]$ is a class of a line in the fiber of the exceptional divisor over $x$, and $i\colon \PP^1_{k(x)}\to Q^{[2]}$ is the embedding. Hence
    \begin{equation}\label{eq:k-coh_hilb}
    \H^{2n-1}(Q^{[2]},\kM_{2n}) = \sum_{\text{$L/k$ finite}} (\pi_{L/k})_* \left(  \kM_1(L) \cdot \H^{2n-1}(Q_L^{[2]},\kM_{2n-1})\right).
    \end{equation}
    Denote by
    \[
    \Gamma:=\mr{Bl}_{\Delta}(Q\times Q)\subseteq Q\times Q^{[2]},
    \]
    the universal subscheme, the projection $f\colon \Gamma \to Q^{[2]}$ is the finite flat degree $2$ quotient morphism, while the projection $g\colon \Gamma \to Q$ is given by the first coordinate. Denote by 
    \[
    E\subseteq \Gamma,\quad E\cong \PP(T_Q),
    \]
    the exceptional divisor and $D:=f(E)\subseteq Q^{[2]}$, $D\cong \PP(T_Q)$, its image under the quotient morphism. These fit into the following commutative diagram.
    \[
    \xymatrix{
    E \ar[r] \ar[d] & \Gamma \ar[d]_f & \Gamma\setminus E \ar[l]\ar[d] & Q\times Q \setminus \Delta \ar[l]_\simeq\ar[d] \\
    D \ar[r] & Q^{[2]} & Q^{[2]}\setminus D \ar[l] & (Q\times Q \setminus \Delta)/(\Z/2) \ar[l]_(0.6)\simeq
    }
    \]
    Here $f|_E\colon E\to D$ is an isomorphism and as Cartier divisors we have
    \begin{equation}\label{eq:Rost_Cartier_ramification}
    f^*D = 2E.
    \end{equation}

    Below we will perform computations with elements of the Gersten complexes for which we use the following notation. For a variety $X$ we put
    \[
    \mc{C}_i(X):=\bigoplus_{x\in X_{(i)}} \K^\M_i(k(x)),\quad \mc{Z}_i(X):=\ker (\partial\colon \mc{C}_i(X)\to \mc{C}_{i-1}(X)),
    \]
    and $\bar{\mc{C}}_i(X)$, $\bar{\mc{Z}}_i(X)$ denote the respective groups in the mod $2$ complex. For $\alpha\in \mc{C}_i(X)$ we denote by $\bar{\alpha}$ its class modulo $2$, for $\alpha\in \mc{Z}_i(X)$ we denote $[\alpha]$ its cohomology class and if $\bar{\alpha}\in \bar{\mc{Z}}_i(X)$ then $[\bar{\alpha}]$ is its cohomology class in mod $2$ Gersten complex. 
 
    Consider an element $[\bar{\alpha}]\in \H^{n-1}(Q,\kM_{n})$ which we want to show to be a sum of pushforwards of decomposable elements. Represent it by
    \[
    \alpha = \sum \{u_i\}_{y_i} \in \mc{C}_1(Q),\quad \bar{\alpha}\in \bar{\mc{Z}}_1(Q),
    \]
    where $y_i\in Q_{(1)}$ are different and $u_i\in k(y_i)^\times$. Without loss of generality we may assume $u_i$ to be transcendental over $k$, otherwise $\{u_i\}_{y_i}$ is obtained as a pushforward $(\pi_{k(u_i)/k})_*$ of a decomposable element and may be subtracted from $\alpha$. We claim that 
    \begin{equation}\label{eq:Rost_pushpull}
    [\bar{\alpha}] \in g_*f^*(\H^{2n-1}(Q^{[2]},\kM_{2n})),
    \end{equation}
    in view of the equality~\eqref{eq:k-coh_hilb} this will give the claim of the theorem.
    
    Denote by $C_i':=\bar{y_i}$ the respective proper curve on $Q$ and $C_i\to C_i'$ its normalization. Since $u_i$ is not algebraic it defines a non-constant morphism $C_i\to \PP^1$. In characteristic $p$ we have $u_i=\tilde{u}_i^{p^{e_i}}$ with $\tilde{u}_i$ such that $d\tilde{u}_i\neq 0$, and 
    \[
    \overline{\{u_i\}}_{y_i}= p^{e_i} \cdot \overline{\{\tilde{u}_i\}}_{y_i}=\overline{\{\tilde{u}_i\}}_{y_i}
    \]
    since $p$ is odd. Changing $u_i$ to $\tilde{u}_i$ we may assume that $u_i\colon C_i\to \PP^1$ are generically \'etale. Consider the cartesian square
    \[
    \xymatrix{
    C_{i,j} \ar[r] \ar[d] & C_j \ar[d]^{u_j} \\
    C_i \ar[r]^{u_i} & \PP^1
    }
    \]
    and denote $C_{i,j}^m$ the normalizations of the irreducible components of $C_{i,j}$. If $i=j$ then there is a diagonal component which we denote $C_{i,i}^\Delta$. Each curve $C_{i,j}^m$ comes with a canonical morphism to $Q\times Q$, and if $C_{i,j}^m\neq C_{i,i}^\Delta$ then its image is not on the diagonal, hence it defines a point $y_{i,j}^m \in (Q\times Q\setminus \Delta)_{(1)}$ giving rise an element
    \[
    \beta':=\sum_{C_{i,j}^m\neq C_{i,i}^\Delta} \{u_i\}_{y_{i,j}^m} \in \mc{C}_1(\Gamma\setminus E)\subseteq \mc{C}_1(\Gamma).
    \]
    The element $\beta'$ is $\Z/2$-invariant, so it descends to $\beta \in \mc{C}_1(Q^{[2]}\setminus D)\subseteq \mc{C}_1(Q^{[2]})$ such that $f^*\beta=\beta'$.

    Since $\bar{\alpha}\in \bar{\mc{Z}}_1(Q)$, we have
    \begin{equation}\label{eq:Rost_divisor_even}
    \sum (\nu_i)_*\mr{div}_0(u_i) - \sum (\nu_i)_*\mr{div}_\infty(u_i) = 2a
    \end{equation}
    for some $a\in \mc{Z}_0(Q)$, here $\nu_i\colon C_i\to \bar{y}_i\to Q$ is the composition of the normalization and the closed embedding. Denote by $d_x^0$ the component of $\sum (\nu_i)_*\mr{div}_0(u_i)$ at $x\in Q_{(0)}$ and similarly $d_x^\infty$ for $\sum (\nu_i)_*\mr{div}_\infty(u_i)$, so we have
    \[
    \sum (\nu_i)_*\mr{div}_0(u_i) = \sum_x d_x^0\cdot [x],\quad \sum (\nu_i)_*\mr{div}_\infty(u_i) = \sum_x d_x^\infty\cdot [x].
    \]
    Put
    \[
    N:=\sum d_x^0 \cdot [k(x):k] = \sum \deg ((\nu_i)_*\mr{div}_0(u_i))
    \]
    and note that $N$ is even by Springer's theorem \cite[Corollary~71.3]{EKM08}. We have
    \[
    g_*(\{u_i\}_{y_{i,j}^m}) = [k(y_{i,j}^m):k(y_i)]\cdot \{u_i\},\quad \sum_{j,m} [k(y_{i,j}^m):k(y_i)] = \sum_{j} [k(y_{j}):k(u_j)] =N,
    \]
    hence viewing $\beta'$ as an element of $\mc{C}_1(\Gamma)$ we obtain
    \begin{equation}\label{eq:Rost_push_Q}
    g_*\beta' =(N-1)\cdot \alpha.
    \end{equation}

    Now we compute $\partial\bar{\beta}$. First we claim that it is supported on $D$. In order to check it we may pass to the algebraic closure of $k$, since $\bar{\mc{C}}_0(Q^{[2]})\to \bar{\mc{C}}_0(Q^{[2]}_{k^{\mr{alg}}})$ is injective. Let $x\in Q^{[2]}(k^{\mr{alg}})$ be the point corresponding to $(x_1,x_2)\in (Q\times Q)(k^{\mr{alg}})$, $x_1\neq x_2$. Then $\partial_x \beta = (d^0_{x_1}d^0_{x_2} - d^\infty_{x_1}d^\infty_{x_2})\cdot [x]$.
    It follows from~\eqref{eq:Rost_divisor_even} that $d^0_{x_i}$ and $d^\infty_{x_i}$ have the same parity, hence $d^0_{x_1}d^0_{x_2} - d^\infty_{x_1}d^\infty_{x_2}$ is even and $\partial_x \bar{\beta} =0$. Thus there exists some $c\in \mc{C}_0(D)$ such that 
    \[
    \partial \bar{\beta} = j_* \bar{c}
    \]
    for the closed embedding $j\colon D\to Q^{[2]}$. Denote by $p\colon D\cong \PP(T_Q)\to Q$ the projection. We claim that $p_*\bar{c} = \bar{a}$ where $a$ was defined in~\eqref{eq:Rost_divisor_even}. Again, this can be checked over the algebraic closure. We have 
    \[
    p_*\bar{c} = \sum_{x\in Q_{(0)}} \tfrac{1}{2}\left(d^0_x(d^0_x-1)- d^\infty_x(d^\infty_x-1)\right)\cdot \overline{[x]}.
    \]
    If $a_x$ is the component of $a$ corresponding to the point $x$, then by~\eqref{eq:Rost_divisor_even} we have $d^0_x-d^\infty_x=2a_x$, thus $\tfrac{1}{2}\left(d^0_x(d^0_x-1)- d^\infty_x(d^\infty_x-1)\right) = 2(a_xd^\infty_x+a_x^2) -a_x$, hence
    \[
    p_*\bar{c} = \bar{a}.
    \]
    
    We have $\partial \alpha =2a$, hence $2[a]=0$ in $\CH_0(Q)$. Since $\CH_0(Q)$ is torsion free \cite[Corollary~71.4]{EKM08}, $[a]=0$. Hence $p_*[\bar{c}]=[\bar{a}]=0$. Pushforward $p_*\colon \CH_0(\PP(T_Q))\to \CH_0(Q)$ is an isomorphism, thus $[\bar{c}]=0$. Then there exists $\gamma\in \mc{C}_1(D)$ such that
    \[
    \partial \bar{\gamma} = \bar{c}.
    \]

    Consider $\beta - j_*\gamma \in \mc{C}_1(Q^{[2]})$. We have $\partial (\bar{\beta}-j_*\bar{\gamma})=j_*\bar{c}-\partial j_*\bar{\gamma}=0$, hence $\bar{\beta} - j_*\bar{\gamma} \in \bar{\mc{Z}}_1(Q^{[2]})$. Furthermore, it follows from~\eqref{eq:Rost_Cartier_ramification} that $f^*(\bar{\beta} - j_*\bar{\gamma})=f^*\bar{\beta}=\bar{\beta}'$. Finally, applying~\eqref{eq:Rost_push_Q} we obtain
    \[
    [\bar{\alpha}] = g_*f^*[\bar{\beta} - j_*\bar{\gamma}] \in g_*f^*(\H^{2n-1}(Q^{[2]},\kM_{2n})).
    \]
    The claim of the theorem follows from~\eqref{eq:k-coh_hilb}.
\end{proof}

\begin{lemma} \label{lem:baby_rost}
    Let $Q$ be a smooth projective quadric of dimension $n\ge 1$ over $k$. Then
    \[
    \Ch^{n-1}(Q) = \Z/2\cdot [h]^{n-1} + \sum_{\substack{\text{$L/k$ finite}\\ i(q_L)\ge 2}} (\pi_{L/k})_*\left( \Z/2\cdot [\ell_L]\right),
    \]
    where $[h]\in \Ch^1(Q)$ is the class of a hyperplane section,  $\pi_{L/k}$ is the projection $Q_L\to Q$ and $[\ell_L]\in \Ch^{n-1}(Q_L)$ is the class of a line.
\end{lemma}
\begin{proof}
    We use the same notation as in the proof of Theorem~\ref{thm:rost}. Consider the universal family
    \[
    Q \xleftarrow{g} \Gamma \xrightarrow{f} Q^{[2]}.
    \]
    
    We claim that $g_*f^*\colon \Ch^{2n-1}(Q^{[2]}) \to \Ch^{n-1}(Q)$ is surjective. Let $y\in Q_{(1)}$ and $C\to \bar{y}\subseteq Q$ be the normalization of the curve $\bar{y}$. Pick a rational function $u\in k(y)$ which is not algebraic over $k$ and which defines a generically \'etale morphism $C\to \PP^1$. Changing $u$ to $u^2$ we may assume that $N:=[k(y):k(u)]$ is even. Consider the cartesian square
    \[
    \xymatrix{
    C\times_{\PP^1}C \ar[r]\ar[d] & C \ar[d]^u\\
    C \ar[r]^u & \PP^1
    }
    \]
    and let $\{y_i\}_{i\in I}$ be the generic points of the irreducible components of $(C\times_{\PP^1}C) \setminus \Delta(C)$. By the construction, $\{y_i\}_{i\in I}$ come equipped with an embedding to $(Q\times Q)\setminus \Delta \subseteq \Gamma$, giving rise to a $1$-cycle $B':=\sum [y_i]$ on $\Gamma$. The set of points $\{y_i\}_{i\in I}\subseteq (Q\times Q)\setminus \Delta$ is $\Z/2$-invariant, thus defines a cycle $B$ on $((Q\times Q)\setminus \Delta)/(\Z/2)\subseteq Q^{[2]}$ such that $f^*[B]=[B']$. We have $g_*[B']=\sum_{i\in I} [k(y_i):k(C)][C]=(N-1)[C]$ in $\CH^{n-1}(Q)$. Since $N$ is even then
    \[
    g_*f^*[B]=[C] \in \Ch^{n-1}(Q).
    \]

    As in the proof of Theorem~\ref{thm:rost} using Lemma~\ref{lem:Hilb_Q} and applying blow-up formula we obtain
    \[
    \Ch^{2n-1}(Q^{[2]})=r^*\Ch^{2n-1}(\Gr(2,V))
    \oplus s_*(h^*\Ch^{2n-3}(\mr{OGr}(2,q))\cdot c_1).
    \]
    
    The first summand is generated by $r^*[\lambda]$, where $\lambda$ is a line on $\Gr(2,V)$ lying in the image of an embedding $\PP(V/L)\subseteq \Gr(2,V)$ for a dimension $1$ subspace $L\leqslant V$. Choose an anisotropic line $L\leqslant V$, a projection from the corresponding point of $\PP(V)$ gives rise to a finite flat degree $2$ morphism $Q\to \PP(V/L)$. By the universal property of the Hilbert scheme we obtain the following commutative diagram with the right-hand square being cartesian.
    \[
    \xymatrix{
    Q \ar[d]_= &Q\ar[d]^{\psi} \ar[l]_= \ar[r]^(0.4){f_\PP} & \PP(V/L) \ar[d]^\varphi\\
    Q & \Gamma \ar[r]^f \ar[l]_g & Q^{[2]} 
    }
    \]
    We have $\varphi(L')=\ell_{L'}\cap Q$ for a point $L'\in \PP(V/L)$, where $\ell_{L'}$ is the line in $\PP(V)$ corresponding to the dimension $2$ preimage of $L'\leqslant V/L$ in $V$, and the intersection is a length $2$ subscheme of $Q$. It is straightforward to check that $r^*[\lambda]=\varphi_*[\ell]$ for a line $\ell\subseteq \PP(V/L)$. Then applying base change we obtain $g_*f^*r^*[\lambda] = f_\PP^*[\ell]$. It remains to notice that $[\ell]=[H]^{n-1}$ for a hyperplane $H\subseteq \PP(V/L)$ and $f^*_\PP([H])=[h]$.
    
    The second summand is generated by the classes $(\pi^{[2]}_{k(x)/k})_*[\ell_x]$, where $\pi^{[2]}_{k(x)/k}\colon Q^{[2]}_{k(x)}\to Q^{[2]}$ is the projection and $[\ell_x]$ is the class of a line in the fiber $\PP^2_{k(x)}$ of the exceptional in $Q^{[2]}$ over a point $x\in \mr{OGr}(2,q)$. Note that $i(q_{k(x)})\ge 2$ since $x$ is a closed point of $\mr{OGr}(2,q)$. We have a commutative diagram
    \[
    \xymatrix{
    Q_{k(x)} \ar[d]_{\pi_{k(x)/k}} & \Gamma_{k(x)} \ar[d]\ar[r]^{f_{k(x)}} \ar[l]_{g_{k(x)}} & Q_{k(x)}^{[2]} \ar[d]^{\pi^{[2]}_{k(x)/k}}\\
    Q & \Gamma \ar[r]^f \ar[l]_g & Q^{[2]}
    }
    \]
    so by base change we have $g_*f^*(\pi^{[2]}_{k(x)/k})_*[\ell_x]=(\pi_{k(x)/k})_*(g_{k(x)})_*f^*_{k(x)}[\ell_x]$. Recall that $(\PP^1)^{[2]}\cong \PP^2$, thus for the embedding of the line $i\colon\PP^1_{k(x)}=l_{k(x)}\to Q_{k(x)}$ corresponding to the rational point $x\in\mr{OGr}(2,q)_{k(x)}$ we have the following commutative diagram with the right-hand square being cartesian.
    \[
    \xymatrix{
    \PP^1_{k(x)} \ar[d]_i & (\PP^1\times \PP^1)_{k(x)} \ar[d] \ar[l]_(0.6){g_{\PP}} \ar[r]^(0.6){f_{\PP}} & \PP^2_{k(x)} \ar[d]^j \\
    Q_{k(x)} & \Gamma_{k(x)} \ar[r]^{f_{k(x)}} \ar[l]_{g_{k(x)}} & Q_{k(x)}^{[2]} 
    }
    \]
    Here $j$ is the embedding of the fiber of the exceptional divisor over $x\in \mr{OGr}(2,q)_{k(x)}$. Then
    \[
    (g_{k(x)})_*f^*_{k(x)}[\ell_x] = i_*(g_{\PP})_*f^*_{\PP} [\ell]
    \]
    for the class $[\ell]$ of a line on $\PP^2_{k(x)}$. We have $f^*_{\PP} [\ell]=[0\times\PP_{k(x)}^1]+[\PP_{k(x)}^1\times 0]$, hence $i_*(g_{\PP})_*f^*_{\PP} [\ell]=i_*1=[\ell_{k(x)}]$ and the claim follows.
\end{proof}

\begin{corollary} \label{cor:rost}
    Let $Q$ be a smooth projective quadric over $k$. Then
    \[
    \H^{n-1}(Q,\kM_n) =  \kM_1(k) \cdot [h]^{n-1}
    + \sum_{\substack{\text{$L/k$ finite}\\ i(q_L)\ge 2}} (\pi_{L/k})_*\left( \kM_1(L)\cdot [\ell_L]\right)
    \]
    where $\pi_{L/k}$ is the projection $Q_L\to Q$, $[h]\in \Ch^1(Q)$ is the class of a hyperplane section and $[\ell_L]\in \Ch^{n-1}(Q_L)$ is the class of a line.
\end{corollary}
\begin{proof}
    Follows as a combination of Theorem~\ref{thm:rost} and Lemma~\ref{lem:baby_rost}.
\end{proof}

\section{Some explicit formulas}\label{sec:explicitformulas}
In this section we give formulas for the completion of the unimodular row $x_\bullet:=(x_1,\hdots,x_{n+1})$ over the ring $R:=k[x_1,\hdots,x_{n+1}]/(q-1)$ for a regular quadratic form $q=\sum_{i=1}^{n+1}a_ix_i^2$ over $k$ satisfying some conditions. 
If $n=1$ or $n=2$ we usually use the notation $x:=x_1$, $y:=x_2$ and $x:=x_1$, $y:=x_2$, $z:=x_3$ respectively.

\begin{example} \label{ex:Revoy}
    Let $q$ be isotropic, then $x_\bullet$ is completable by {\cite[2.1]{Rev79}} as follows. Choose a basis in such a way that $q$ becomes $\tilde{q}=X_1X_2+\sum_{i=3}^{n+1} b_iX_i^2$ (see e.g. \cite[Theorem~I.3.4]{Lam05}), denote
    \begin{gather*}
    \lambda:=1-X_2,\quad \mu:=(1-X_2)X_2-1,\quad U:=(X_3,X_4,\hdots,X_{n+1}),\\
    V:=(b_3X_3,b_4X_4,\hdots,b_{n+1}X_{n+1}),
    \end{gather*}
    and consider the matrix
    \[
        \tilde{M}:=
        \begin{pmatrix}
            X_1 & X_2 & U \\
            -1 & \lambda X_2 & 0\\
            -\lambda V^T & \lambda\mu V^T & 1_{n-1}
        \end{pmatrix}.
    \]
    A straightforward computation shows that $\det \tilde{M}= (1-X_2)\tilde{q}+X_2 = 1$. In order to get a completion for the original $x_\bullet$ over $R$ one has to transform this matrix to the original basis.

    In dimension $n=2$ one can also proceed as follows. Since $q$ is isotropic, $a_1t_1^2+a_2t_2^2+a_3t_3^2=0$ for some nontrivial $t_1,t_2,t_3\in k$, and without loss of generality we may assume $t_1\in k^\times $. Then for
    \[
        X_1:=\frac{1}{2t_1}x_1-\frac{a_2t_2}{2a_1t_1^2}x_2-\frac{a_3t_3}{2a_1t_1^2}x_3,
        \,
        X_2:=2a_1t_1x_1+2a_2t_2x_2+2a_3t_3x_3,
        \,
        X_3:=\frac{t_3}{2a_1t_1}x_2-\frac{t_2}{2a_1t_1}x_3,
    \]
    one has
    \[
        \det \begin{pmatrix}
            x_1 & x_2 & x_3\\
            t_1 & t_2 & t_3\\
            -\frac{a_2a_3}{t_1}X_3 & \frac{a_2a_3t_2}{t_1^2} X_3 + \frac{2a_3t_3}{t_1}X_1 & \frac{a_2a_3t_3}{t_1^2} X_3- \frac{2a_2t_2}{t_1}X_1
        \end{pmatrix}
        = 1.
    \]
    This matrix is obtained from
    \[
        \begin{pmatrix}
            X_1 & X_2 & X_3\\
            1 & 0 & 0\\
           0 & bX_3 & -X_1
        \end{pmatrix},\quad  b:=-4a_1a_2a_3 ,
    \]
    by the change of basis transforming $X_1X_2+bX_3^2$ to $q$.
 \end{example}

\begin{example}
    If $n=1$ then $x_\bullet=(x,y)\in R^2$ is completable as 
    \[
    \begin{pmatrix}
        x & y \\
        -a_2y & a_1x
    \end{pmatrix}.
    \]
\end{example}

\begin{example} \label{ex:dim2_1}
    The proof of Theorem~\ref{thm:hedge2} provides an explicit recipe for the completion of the unimodular row $x_\bullet=(x,y,z)$ over $R$ with the entries of the matrix of degree $\le 1$ in $x$, $y$ and $z$. The particular answer depends heavily on the choice of the quaternion basis $\{1,u,v,uv\}$, which in turn may be made explicit based on $(x_0,x_1,x_2,x_3)$, but the formulas are cumbersome. For example, for $k=\mathbb{Q}[\omega]$, $\omega:=\sqrt{-7}$, and $q=x^2+y^2+z^2$ choosing 
    \[
    u:=i+2j-\omega ij,\quad v:= -\frac{1}{2}i+2j-\frac{\omega}{2}ij,\quad uv=\omega i + \omega j +3ij,
    \]
    and applying some elementary transformations on the last two rows of the final matrix we obtain
    \[
        M':=
        \begin{pmatrix}
        x & y & z \\
        2-\omega (x-y)-2z & 10-\omega(x-y)-2z & -4\omega-x+5y\\
        \omega(x-y)+5z & 6+4\omega x+11z & -2\omega+8x-y-4\omega z
        \end{pmatrix}
    \]
    with $\det M'=-40$. Compare this to \cite[Theorem~1.3]{Mue26}.
\end{example}

\begin{example} \label{ex:dim2_2}
    Let $q=-x^2+ay^2+az^2$, then over $R$ we have
    \[
        (x,y,z)
        \begin{pmatrix}
            2 & 0 & 0\\
            ay & 1 & 0\\
            az & 0 & 1
        \end{pmatrix}
        =
        ((x+1)^2,y,z).
    \]
    The row $(x+1,y,z)$ is unimodular with
    \[
        \tfrac{1}{2}(1-x)\cdot (x+1)+\tfrac{1}{2}ay\cdot y+ \tfrac{1}{2} az \cdot z = 1,
    \]
    thus $((x+1)^2,y,z)$ may be completed to an invertible matrix as
    \[
        \begin{pmatrix}
                (x+1)^2 & y & z \\
                -y-a(x+1)z & \frac{a^2z^2}{4} & \frac{2-2x-a^2yz}{4} \\
                -z+a(x+1)y & \frac{2x-2-a^2yz}{4} & \frac{a^2y^2}{4}
        \end{pmatrix}
    \]
    (see \cite[III.(4.4)]{Lam06} for the particular formula that we used). Transforming back to the original row we obtain an invertible matrix
    \[
        \begin{pmatrix}
            x & y & z \\
            -\frac{2y+a(x+3)z}{4} & \frac{a^2z^2}{4} & \frac{2-2x-a^2yz}{4} \\[1mm]
            \frac{-2z+a(x+3)y}{4} & \frac{2x-2-a^2yz}{4} & \frac{a^2y^2}{4}
        \end{pmatrix}.
    \]
\end{example}

\begin{example} \label{ex:dim2_3}
    Let $q=ax^2+by^2-abz^2$. Consider the quaternion algebra $\mr{H}:=(a,b)_k$ and put
    \[
        \mr{H}(R):= R \otimes_k \mr{H},\quad \mr{H}^0(R):= R \otimes_k \mr{H}^0 = Ri\oplus Rj\oplus Rij, \quad \rho:=xi+yj+zij\in \mr{H}^0(R).
    \] 
    Consider the non-degenerate symmetric bilinear form
    \[
    \mr{H}^0(R)\times\mr{H}^0(R)\to R,\quad (h_1,h_2)\mapsto \langle h_1,h_2\rangle := \operatorname{Trd}(h_1h_2).
    \]
    Then decomposition~\eqref{eq:tangent_stably_trivial} becomes
    \[
    \mr{H}^0(R) = \T(R)\oplus R\rho,\quad \T(R):=(R\rho)^\perp=\{h\in \mr{H}^0(R)\,|\,\operatorname{Trd}(\rho h)=0\}
    \]
    and $\langle-,-\rangle$ restricts to a non-degenerate symmetric bilinear form on $\T(R)$.
    Note that $\rho^2=1$, $\rho \T(R) = \T(R)$ and $\rho h = -h\rho$ for every $h\in\T(R)$, hence
    \[
    \T(R) = \T(R)_+\oplus \T(R)_-, \quad \T(R)_+:= \frac{1 + \rho}{2} \T(R),\quad \T(R)_-:= \frac{  1 - \rho}{2}\T(R)
    \]
    with both modules $\T(R)_\pm$ being rank $1$ projective modules over $R$ isotropic with respect to $\langle-,-\rangle$. Since the form is non-degenerate then it induces an isomorphism $\T(R)_+\cong \T(R)_-^\vee$. Consider the elements
    \[
    e^+_1:= \frac{1 + \rho}{2}(i-ax\rho),\quad e^+_2:= \frac{1 + \rho}{2}(j-by\rho),\quad e^-_1:= \frac{1 - \rho}{2}(i-ax\rho),\quad e^-_2:= \frac{1 - \rho}{2}(j-by\rho).
    \]
    These are $i$ and $j$ projected first to $\T$ and then further to $\T_+$ and $\T_-$ respectively. The Gram matrix of these vectors is the following one.
    \[
    G(e^+_1,e^+_2,e^-_1,e^-_2)=\frac{1}{2}\begin{pmatrix}
        0 & 0 & a(1-ax^2) & -ab(xy+z) \\
        0 & 0 & ab(z-xy) & b(1-by^2) \\
        a(1-ax^2) & ab(z-xy) & 0 & 0 \\
        -ab(xy+z) & b(1-by^2) & 0 & 0
        \end{pmatrix}.
    \]
    We want to find a basis for $\T(R)$ as $v:=e_1^++\lambda_1e_1^-+\mu_1e_2^-$, $w:=e_2^++\lambda_2e_1^-+\mu_2e_2^-$, and for this it is necessary and sufficient to find $\lambda_1,\mu_1,\lambda_2,\mu_2\in R$ such that the Gram matrix for $v,w$ is invertible. A straightforward computation shows that
    \[
    \det G(v,w)=-\frac{a^2b^2}{4}(\mu_1(x^2-bz^2)-\lambda_2(y^2-az^2)+(\mu_2-\lambda_1)xy+(\lambda_1+\mu_2)z)^2.
    \]
    We would like to make the formula in the brackets a nonzero scalar multiple of $-ax^2-by^2+abz^2$, and this can be achieved with
    \[
    \lambda_1=\mu_2=abz,\quad \mu_1=2a,\quad \lambda_2=-2b.
    \]
    For this choice of parameters we have $\det G(v,w)=-a^2b^2$, hence $v,w$ form a basis of $\T(R)$ yielding an explicit completion of $(x,y,z)$
    \[
    A:=\begin{pmatrix}
    x & y & z \\
    l_1 & l_2 & l_3 \\
    m_1 & m_2 & m_3
    \end{pmatrix},
    \]
    of determinant $\det A=4$ with 
    \[
    \begin{aligned}
    l_1&=(1+abz)(1-ax^2)-2ab(xy+z),& m_1&=b(z-xy)-2b(1-ax^2)-ab^2z(xy+z)\\
    l_2&=-a(xy+z)+a^2bz(z-xy)+2a(1-by^2),& m_2&=(1+abz)(1-by^2)-2ab(z-xy),\\
    l_3&=-(axz+y)+abz(y-axz)-2a(x+byz), &  m_3&=x-byz-2b(y-axz)-abz(x+byz).
    \end{aligned}
    \]
\end{example}

\begin{example}\label{ex:dim3}
    Let $q=\sum_{i=1}^4 a_ix_i^2$ and suppose that $a_1a_2a_3a_4=d^2$ for some $d\in k^\times$. Then $x_\bullet$ is completable as
    \[
        M=
        \begin{pmatrix}
        x_1 & x_2 & x_3 & x_4 \\
        -\frac{a_2}{a_1}x_2
            & x_1
            & -\frac{d}{a_1a_3}x_4
            & \frac{d}{a_1a_4}x_3 \\
        -a_3x_3
            & \frac{d}{a_2}x_4
            & a_1x_1
            & -\frac{d}{a_4}x_2 \\
        -a_4x_4
            & -\frac{d}{a_2}x_3
            & \frac{d}{a_3}x_2
            & a_1x_1
        \end{pmatrix}.
    \]
    Here $\det M = q^2=1$ and the matrix is a slight modification of the one from the proof of \cite[Theorem~6.11]{Fas11a}, where the second, third and fourth rows are the natural choices of vectors orthogonal to $(x_1,x_2,x_3,x_4)$.
    
    In the particular case of a Pfister form $q=x_1^2-ax_2^2-bx_3^2+abx_4^2$ the matrix looks as follows.
    \[
        M=
        \begin{pmatrix}
            x_1 & x_2 & x_3 & x_4\\
            a x_2 & x_1 & a x_4 & x_3\\
            b x_3 & -b x_4 & x_1 & -x_2\\
            -ab x_4 & b x_3 & -a x_2 & x_1
        \end{pmatrix}.
    \]
    See also \cite[2.2]{Rev79}.
\end{example}

\end{document}